\documentclass{amsproc} %article [reqno]<-equation numbers on right of page
\usepackage{amsmath,mathrsfs,amsthm,amsfonts}%,graphicx,amscd,amssymb,bbm}
\usepackage[margin=1in]{geometry}
\usepackage{enumitem}
\usepackage[figuresright]{rotating} 
\usepackage{tikz,siunitx}
\usepackage{tikzit}
\tikzstyle{skein}=[-, draw=red, line width=2]
\tikzstyle{base}=[-, draw=black, line width=2, tikzit draw=blue]
\tikzstyle{arrow1}=[draw={rgb,255: red,128; green,128; blue,128}, ->]
\tikzstyle{line}=[-, draw={rgb,255: red,128; green,128; blue,128}]
\tikzstyle{arrow2}=[draw={rgb,255: red,128; green,128; blue,128}, <-]

\tikzstyle{reddot1}=[draw=red, shape=circle, fill=red, inner sep=0pt, minimum size=3pt]

\tikzstyle{fill}=[-]
\tikzstyle{redline}=[-, fill=none, draw=red]
\tikzstyle{blueline}=[-, draw=blue]
\tikzstyle{greenline}=[-, draw={rgb,255: red,25; green,255; blue,0}]
\tikzstyle{greenline1}=[-, draw={rgb,255: red,0; green,128; blue,36}]
\tikzstyle{arrow1}=[->, draw={rgb,255: red,11; green,162; blue,0}]
\tikzstyle{dashed1}=[-, dashed]
\tikzstyle{arrow2}=[->]
\tikzstyle{arrow3}=[<-]
\tikzstyle{orange1}=[-, draw={rgb,255: red,255; green,128; blue,0}]
\tikzstyle{arrow4}=[<-, draw={rgb,255: red,28; green,118; blue,3}]
\tikzstyle{purple1}=[-, draw={rgb,255: red,217; green,29; blue,208}]
\tikzstyle{brown1}=[-, draw={rgb,255: red,162; green,108; blue,54}]

\usepackage{array}
\usepackage{algorithm}
\usepackage[noend]{algpseudocode}
\usetikzlibrary{matrix,positioning,calc,arrows,decorations.pathmorphing}
\usepackage{wrapfig}
\usepackage{graphicx}
\usepackage{caption} 
\usepackage{multirow}
\theoremstyle{plain}
\newtheorem{Th}{Theorem}[section]
\newtheorem{Thm}[Th]{Theorem}

\newtheorem{Lem}[Th]{Lemma}
\theoremstyle{definition}
\newtheorem{Rem}[Th]{Remark}
\newtheorem{Def}[Th]{Definition}

\newtheorem{Prop}[Th]{Proposition}

\allowdisplaybreaks[4]
\pgfdeclarelayer{edgelayer}
\pgfdeclarelayer{nodelayer}
\pgfsetlayers{edgelayer,nodelayer,main}
\tikzstyle{none}=[inner sep=0pt]
\begin{document}
\title{Algorithms for Skein Manipulation in a Genus-2 Handlebody}

%Authors are listed in alphabetical order by last name.

\author{Rachel Kinard}
\address{Air Force Research Laboratory, Sensors Directorate (AFRL/RYAT), 2241 Avionics Cir, WP AFB, OH 45433}
\email{rachel.kinard@us.af.mil}

\author{R\u{a}zvan Gelca}
\address{Texas Tech University}
\email{rgelca@gmail.com}

\author{Paul T. Schrader}
\address{Air Force Research Laboratory, Information Directorate (AFRL/RIGC), 525 Brooks Road, Rome, NY 13441}
\email{paul.schrader.1@us.af.mil}

%\thanks will become a 1st page footnote.
%\thanks{}

%General info
\maketitle

\begin{abstract} 
	We present a series of algorithms for skein manipulation in a genus-2 handlebody, implementing a novel strand sorting method to reduce any skein to a skein in a 2-punctured disk. This reduction guarantees resolution as a linear combination of basis elements of the Kauffman Bracket Skein Module. Manually, these skein manipulations prove to be computationally intensive due to the inherent exponential nature of skein relations (i.e., a skein diagram with $n$ crossings yields $2^n$ new skein diagrams, each in $\mathbb{C}[t,t^{-1}]$, the Laurent polynomials with complex coefficients). Thus, as the number of crossings in a skein diagram increases, manual computations become intractable and automation desirable. We enable the automation of all skein computations in the genus-2 handlebody by first converting the skein diagram into an equivalent array, reducing the task of performing skein computations to that of implementing array operators, and then proving that we can always recover the resulting complex Laurent polynomial.
\end{abstract}

%%%%%%%%%%%%%%%%%%%%%%%%%%%%%%%%%%%%%%%%%%%%%%%%%%%%%%%%%%
%%%%%%%%%%%%%%%%%%%%%%%%%%%%%%%%%%%%%%%
\section{Introduction} \label{sec1}
%%%%%%%%%%%%%%%%%%%%%%%%%%%%%
\subsection{Motivation}
\noindent Traditional skein manipulation (the art of producing algebraic representations from a diagrammatic array) is a delicate combination of creativity, care, and complexity employing very carefully executed manual calculations. Knots with only a few crossings (e.g., the trefoil) still require significant computational effort. Even an expert in these computations must be vigilant; the introduction of so much as a minor error or a faulty assumption risks propagation throughout the entire computation, producing incorrect resulting polynomial representations. Skein computations are becoming increasingly important to scientific application in knot theoretical studies, quantum field theory, and quantum information science, arising from the study of knot and link invariants. In these applications, the exponential nature of skein computations and the scale at which these computations must be performed makes automation desirable, producing consistent, reliable results and refocusing the mathematician's effort on more productive theoretical pursuits. Within skein theory exist surprisingly simple patterns and remarkable relationships awaiting discovery; it often happens that a skein computation involving numerous terms simplifies to a relation with very few terms.

%In this case, we do not automate ourselves out of a job; we simply automate ourselves away from computations and into a more desirable job of theoretical progress.

%\noindent The manual calculations used in traditional skein manipulations (the art of producing algebraic representations from a diagrammatic input) are a delicate combination of creativity, care, and complexity. Knots with only a few crossings (e.g., the trefoil) still require significant computational effort. Even an expert manipulator must be vigilant; the introduction of even a minor error risks propagation throughout the computation resulting in catastrophic failure. Additionally, skein manipulation is becoming increasingly important to scientific application in quantum information science and knot theoretical studies. Often in these cases turning to automation and algorithm makes for dependable results and refocuses the mathematician's efforts to more productive pursuits in the theoretical. Here is such a case.

%\subsection{History}
%Content for Paper 1 Background
\subsection{History}
\noindent Historically, skein theory began with John Conway's discovery of a skein relation for the Alexander polynomial, and was formalized by Vaughn Jones' knot polynomial in 1984 \cite{jones1985}. Unlike the Alexander polynomial, the Jones polynomial could only be defined by skein relations. Following Jones, Louis Kauffman defined a relation and link invariant, the Kauffman Bracket, also computed via skein relations \cite{kauffman1987}. Edward Witten \cite{witten1989} gave an intrinsic definition of the Jones polynomial using quantum field theory based on the Chern-Simons functional (1989), and this theory was rigorously developed by Nikolai Reshetikhin and Vladmir Turaev \cite{reshetikhin1991}, \cite{turaev2016}. (We note that there is a parallel theory constructed by Christian Blanchet, Nathan Habegger, Gregor Masbaum, and Pierre Vogel which makes use of the Kauffman bracket in \cite{blanchet1992}; see also \cite{lickorish1993}). In particluar, we replicate the skein computations in the genus 2 handlebody of J. Przytycki \cite{przytycki1991}. Skein computations in the genus-2 handlebody are well understood, however, the ideal way to demonstrate the computational dilemmas arising in even these simple skein manipulations is by providing some rudimentary examples.

\subsection{Our Setting} \label{sec_setting}
\begin{center}
	\begin{figure}[!ht]
%\begin{wrapfigure}{L}[0pt]{6.3cm}
        \input{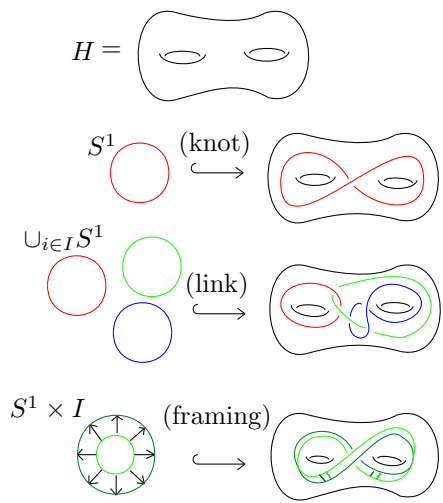}
        \caption{Skeins embedded in a genus-2 handlebody, H}
        \label{fig1}
%\end{wrapfigure}
	\end{figure}
\end{center}

 To begin, let $H$ be a genus-2 handlebody. A \emph{knot} in $H$ is an embedding of a circle $\mathbb{S}^{1} \hookrightarrow H$. A \emph{link} is an embedding of several disjoint copies of $\mathbb{S}^{1}$, ${\bigsqcup}_{i\in I}\mathbb{S}^{1} \hookrightarrow H$, where the index set $I$ numbers the components of the link. We refer to these embeddings as ``skeins''. Similar to skeins of yarn, they are equivalent up to \emph{ambient isotopy} (strands can be manipulated isotopically by combinations of Reidemeister moves without affecting the result). Figure \ref{fig1} shows these embeddings in the genus-2 handlebody. 
 
 Recall the three Reidemeister moves shown in Figure \ref{fig2}. Application of Reidemeister moves can affect the number of crossings we see in a given projection of the skein system. However, the projection itself also affects the seen number of crossings and any computation or series of computations we perform has to be consistent regardless of the projection we choose. 
 
 A \emph{framed knot} in $H$ is the embedding of an annulus, $S^{1} \times I \hookrightarrow H$, where $I$ is the unit interval. It is customary to represent a framed knot by the same diagram as a unframed knot with the convention that the framing is parallel to the plane of the diagram (the "blackboard" framing). We manipulate framed knots in the same manner (by Reidemeister moves) except that we must take the framing into account whenever we twist/untwist any portion of the skein (R1).

Skein computations are traditionally performed diagrammatically. A skein relation, shown in Figure \ref{fig3}, resolves an apparent crossing in a skein as a diagrammatic combination of two other skeins. Resolving all crossings in a skein produces a Laurent polynomial with skein terms. This requires $2^{n}$ computations, where $n$ is the number of crossings in the projection of the original skein. Thus, the resulting Laurent polynomial with coefficients in $\mathbb{C}$, or $\mathbb{C}[t,t^{-1}]$, has $2^{n}$ skein terms. For some of the resulting skeins, we obtain curves which are basis elements of the Kauffman Bracket Skein Module of the genus-2 handlebody (Przytycki); these fundamental curves are shown in Figure \ref{fig4}: $x,y,z$. However, in many cases, resolving all crossings in a skein produces a curve which has no self-intersections, but is not a fundamental curve, such as the skein terms at the botton of Figure \ref{fig3}. For such curves, how do we proceed? In section \ref{sec_method}, we propose a method for resolving all skeins in the genus-2 handlebody.

\begin{center}
	\begin{figure}[!ht]
		\begin{tikzpicture}
	\begin{pgfonlayer}{nodelayer}
		\node [style=none] (0) at (-8.75, 1.5) {};
		\node [style=none] (1) at (-8.75, 0.5) {};
		\node [style=none] (2) at (-7.75, 0.5) {};
		\node [style=none] (3) at (-8.5, 0.2) {};
		\node [style=none] (4) at (-8.5, 0) {};
		\node [style=none] (5) at (-8.5, -1) {};
		\node [style=none] (6) at (-7.5, 0.25) {};
		\node [style=none] (7) at (-7, 0.25) {};
		\node [style=none] (8) at (-7.025, 1.5) {};
		\node [style=none] (9) at (-6.525, 0.125) {};
		\node [style=none] (10) at (-7.025, -1) {};
		\node [style=none] (11) at (-3.75, 1.5) {};
		\node [style=none] (12) at (-3.75, -1) {};
		\node [style=none] (13) at (-2.95, 1.45) {};
		\node [style=none] (14) at (-3.25, 1.05) {};
		\node [style=none] (15) at (-3.375, 0.875) {};
		\node [style=none] (16) at (-3.5, -0.5) {};
		\node [style=none] (17) at (-3.325, -0.675) {};
		\node [style=none] (18) at (-3, -1) {};
		\node [style=none] (19) at (-2.75, 0.25) {};
		\node [style=none] (20) at (-2.25, 0.25) {};
		\node [style=none] (21) at (-2.25, 1.5) {};
		\node [style=none] (22) at (-2.25, -1) {};
		\node [style=none] (23) at (-1.5, 1.5) {};
		\node [style=none] (24) at (-1.5, -1) {};
		\node [style=none] (25) at (1.75, 1.5) {};
		\node [style=none] (26) at (3, -1) {};
		\node [style=none] (27) at (3, 1.5) {};
		\node [style=none] (28) at (2.425, 0.4) {};
		\node [style=none] (29) at (2.325, 0.175) {};
		\node [style=none] (30) at (1.75, -1) {};
		\node [style=none] (31) at (3.75, 1.5) {};
		\node [style=none] (32) at (5, -1) {};
		\node [style=none] (33) at (5, 1.5) {};
		\node [style=none] (34) at (4.425, 0.4) {};
		\node [style=none] (35) at (4.325, 0.175) {};
		\node [style=none] (36) at (3.75, -1) {};
		\node [style=none] (37) at (3, 0.25) {};
		\node [style=none] (38) at (3.75, 0.25) {};
		\node [style=none] (39) at (2.25, 1.5) {};
		\node [style=none] (40) at (1.975, 1.2) {};
		\node [style=none] (41) at (1.85, 1.1) {};
		\node [style=none] (42) at (1.825, -0.575) {};
		\node [style=none] (43) at (1.975, -0.725) {};
		\node [style=none] (44) at (2.25, -1) {};
		\node [style=none] (45) at (4.5, 1.5) {};
		\node [style=none] (46) at (4.75, 1.275) {};
		\node [style=none] (47) at (4.9, 1.1) {};
		\node [style=none] (48) at (4.9, -0.575) {};
		\node [style=none] (49) at (4.75, -0.75) {};
		\node [style=none] (50) at (4.5, -1) {};
		\node [style=none] (51) at (8, 1.5) {};
		\node [style=none] (52) at (8, 0.25) {};
		\node [style=none] (53) at (9.25, 1.5) {};
		\node [style=none] (54) at (9.25, 0.25) {};
		\node [style=none] (55) at (8, -1) {};
		\node [style=none] (56) at (9.25, -1) {};
		\node [style=none] (57) at (8.25, 1.25) {};
		\node [style=none] (58) at (9.1, 1.25) {};
		\node [style=none] (59) at (8.25, 0.85) {};
		\node [style=none] (60) at (9, 0.85) {};
		\node [style=none] (61) at (8.25, 0.425) {};
		\node [style=none] (62) at (9, 0.425) {};
		\node [style=none] (63) at (8.25, 0.05) {};
		\node [style=none] (64) at (9.125, 0.05) {};
		\node [style=none] (65) at (8.25, -0.3) {};
		\node [style=none] (66) at (9.125, -0.3) {};
		\node [style=none] (67) at (8.25, -1) {};
		\node [style=none] (68) at (9, -1) {};
		\node [style=none] (69) at (8.25, -0.65) {};
		\node [style=none] (70) at (9.075, -0.65) {};
		\node [style=none] (71) at (-7.75, -1.6) {\scriptsize{``twist/untwist''}};
		\node [style=none] (72) at (-7.75, -2.175) {\scriptsize{R1}};
		\node [style=none] (73) at (-2.5, -1.55) {\scriptsize{``cross/uncross''}};
		\node [style=none] (74) at (-2.5, -2.125) {\scriptsize{R2}};
		\node [style=none] (75) at (3.25, -1.575) {\scriptsize{``pass/unpass''}};
		\node [style=none] (76) at (3.25, -2.15) {\scriptsize{R3}};
		\node [style=none] (77) at (8.675, -1.65) {\scriptsize{isotopy}};
	\end{pgfonlayer}
	\begin{pgfonlayer}{edgelayer}
		\draw (0.center) to (1.center);
		\draw [bend left=90, looseness=1.50] (2.center) to (1.center);
		\draw [in=90, out=90, looseness=2.00] (2.center) to (3.center);
		\draw (5.center) to (4.center);
		\draw [style=arrow1] (6.center) to (7.center);
		\draw [in=75, out=-90] (8.center) to (9.center);
		\draw [in=90, out=-120] (9.center) to (10.center);
		\draw [bend left=45, looseness=1.25] (11.center) to (12.center);
		\draw (13.center) to (14.center);
		\draw [in=-135, out=135] (16.center) to (15.center);
		\draw (18.center) to (17.center);
		\draw [style=arrow1] (19.center) to (20.center);
		\draw [in=75, out=-75] (21.center) to (22.center);
		\draw [in=105, out=-105] (23.center) to (24.center);
		\draw (25.center) to (26.center);
		\draw (27.center) to (28.center);
		\draw (30.center) to (29.center);
		\draw (31.center) to (32.center);
		\draw (33.center) to (34.center);
		\draw (36.center) to (35.center);
		\draw [style=arrow1] (37.center) to (38.center);
		\draw (39.center) to (40.center);
		\draw (43.center) to (44.center);
		\draw (45.center) to (46.center);
		\draw (49.center) to (50.center);
		\draw [in=135, out=-135] (41.center) to (42.center);
		\draw [in=45, out=-60, looseness=0.75] (47.center) to (48.center);
		\draw [in=60, out=-135, looseness=1.50] (51.center) to (52.center);
		\draw [in=150, out=-15] (53.center) to (54.center);
		\draw [in=-135, out=45, looseness=1.25] (55.center) to (52.center);
		\draw [in=-45, out=120] (56.center) to (54.center);
		\draw [style=arrow1] (57.center) to (58.center);
		\draw [style=arrow1] (59.center) to (60.center);
		\draw [style=arrow1] (61.center) to (62.center);
		\draw [style=arrow1] (63.center) to (64.center);
		\draw [style=arrow1] (65.center) to (66.center);
		\draw [style=arrow1] (67.center) to (68.center);
		\draw [style=arrow1] (69.center) to (70.center);
	\end{pgfonlayer}
\end{tikzpicture}
		\caption{Reidemeister Moves: $R1$, $R2$, and $R3$. Skeins are considered up to isotopy.}
            \label{fig2}
	\end{figure}
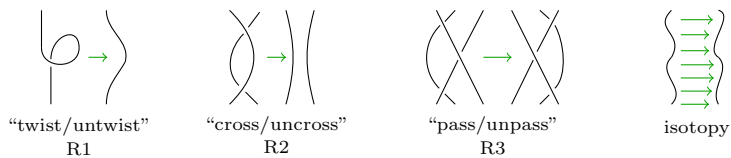
\end{center}

\begin{center}
	\begin{figure}[!ht]
		\input{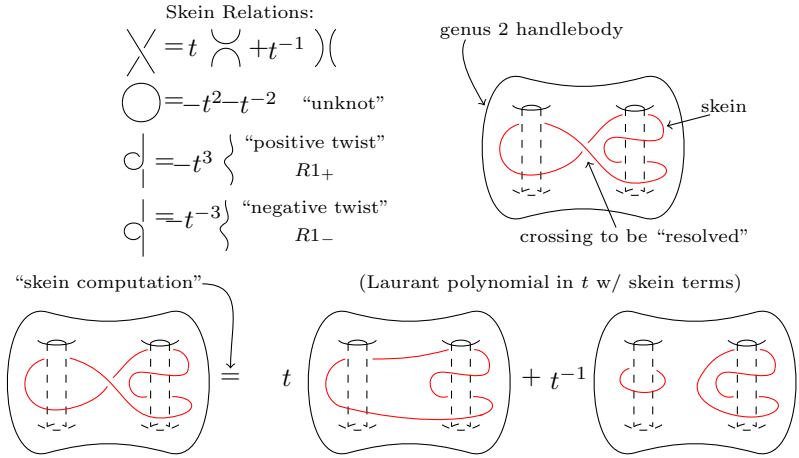}
		\caption{Skein Relations (framed) and Skein Computations; resolving a crossing produces two skein terms, each with a coefficient in $\mathbb{C}[t,t^{-1}]$.}
		\label{fig3}
	\end{figure}
\end{center}

\begin{wrapfigure}{L}[0pt]{5cm}
    \begin{tikzpicture}
	\begin{pgfonlayer}{nodelayer}
		\node [style=none] (0) at (-1.75, 2.25) {};
		\node [style=none] (1) at (1.5, 2.25) {};
		\node [style=none] (2) at (-1.75, -2.25) {};
		\node [style=none] (3) at (1.5, -2.25) {};
		\node [style=none] (4) at (-2, 1.5) {};
		\node [style=none] (5) at (-1, 1.5) {};
		\node [style=none] (6) at (-1.75, 1.425) {};
		\node [style=none] (7) at (-1.25, 1.425) {};
		\node [style=none] (8) at (0.75, 1.5) {};
		\node [style=none] (9) at (1.75, 1.5) {};
		\node [style=none] (10) at (1, 1.425) {};
		\node [style=none] (11) at (1.5, 1.425) {};
		\node [style=none] (12) at (-2, -1.5) {};
		\node [style=none] (13) at (-1, -1.5) {};
		\node [style=none] (14) at (-1.75, -1.55) {};
		\node [style=none] (15) at (-1.25, -1.55) {};
		\node [style=none] (16) at (0.75, -1.5) {};
		\node [style=none] (17) at (1.75, -1.5) {};
		\node [style=none] (18) at (1, -1.55) {};
		\node [style=none] (19) at (1.5, -1.55) {};
		\node [style=none] (20) at (-1.825, 1.025) {};
		\node [style=none] (21) at (-1.15, 1) {};
		\node [style=none] (23) at (-1.15, 0.125) {};
		\node [style=none] (24) at (0.875, -0.725) {};
		\node [style=none] (25) at (1.625, 0.25) {};
		\node [style=none] (26) at (-1.875, 0.1) {};
		\node [style=none] (27) at (1.6, 1.25) {};
		\node [style=none] (29) at (0.825, 0.225) {};
		\node [style=none] (31) at (1.625, -0.725) {};
		\node [style=none] (32) at (-0.75, 1) {};
		\node [style=none] (33) at (3, 2) {};
		\node [style=none] (34) at (2, 0) {};
		\node [style=none] (35) at (3.25, 0.5) {};
		\node [style=none] (36) at (2, -1) {};
		\node [style=none] (37) at (3.5, -1.25) {};
		\node [style=none] (38) at (3.225, 2.05) {$x$};
		\node [style=none] (39) at (3.525, 0.625) {$y$};
		\node [style=none] (40) at (3.8, -1.25) {$z$};
	\end{pgfonlayer}
	\begin{pgfonlayer}{edgelayer}
		\draw [in=-165, out=-15] (0.center) to (1.center);
		\draw [in=-15, out=15] (1.center) to (3.center);
		\draw [in=15, out=165] (3.center) to (2.center);
		\draw [in=165, out=-165] (2.center) to (0.center);
		\draw [bend right=60, looseness=0.50] (4.center) to (5.center);
		\draw [bend left=45, looseness=0.75] (6.center) to (7.center);
		\draw [bend right=45, looseness=0.75] (8.center) to (9.center);
		\draw [bend left=45, looseness=0.75] (10.center) to (11.center);
		\draw [style=dashed1] (6.center) to (14.center);
		\draw [style=dashed1] (7.center) to (15.center);
		\draw [style=dashed1] (10.center) to (18.center);
		\draw [style=dashed1] (11.center) to (19.center);
		\draw [style=dashed1, bend left] (14.center) to (15.center);
		\draw [style=dashed1, bend left] (18.center) to (19.center);
		\draw [style=dashed1, bend right] (16.center) to (17.center);
		\draw [style=dashed1, bend right] (12.center) to (13.center);
		\draw [style=orange1, in=-15, out=-165, looseness=4.50] (20.center) to (21.center);
		\draw [style=blueline, in=-45, out=-150, looseness=0.75] (26.center) to (25.center);
		\draw [style=blueline, in=-165, out=0] (23.center) to (29.center);
		\draw [style=greenline, in=-150, out=-30, looseness=3.75] (31.center) to (24.center);
		\draw [style=arrow4, in=-150, out=0] (32.center) to (33.center);
		\draw [style=arrow4, in=-135, out=15] (34.center) to (35.center);
		\draw [style=arrow4, bend right=15] (36.center) to (37.center);
	\end{pgfonlayer}
\end{tikzpicture}
    \caption{Fundamental Curves}
    \label{fig4}
\end{wrapfigure}
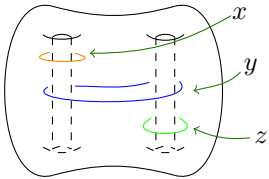
\begin{center}
	\begin{figure}[!ht]
		\input{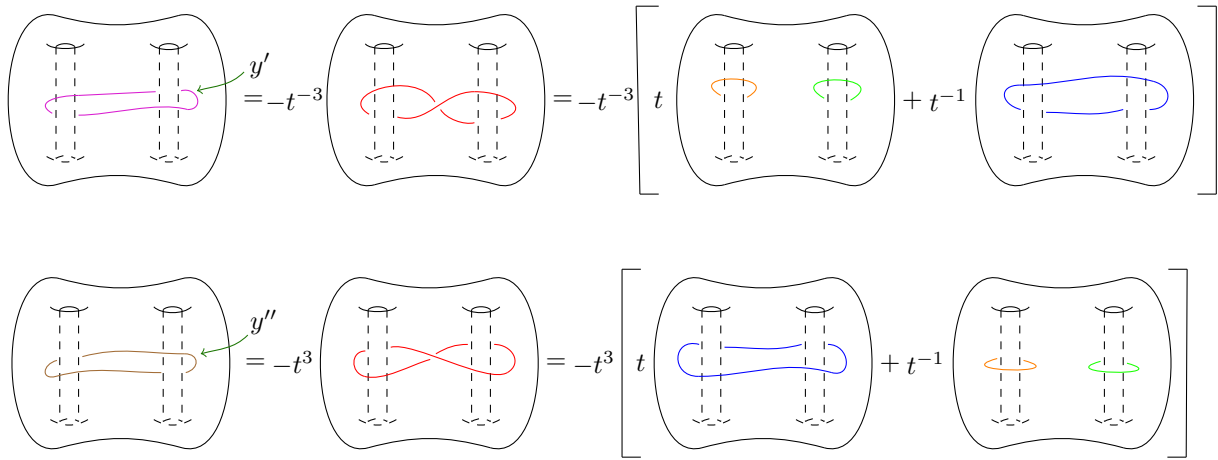}
		\caption{Auxilary curves: $y'=-t^{-2}xz-t^{-4}y$ and $y''=-t^{4}y-t^{2}xz$. These curves occur frequently, and are simple combinations of fundamental curves.}
		\label{fig5}
	\end{figure}
\end{center}

As a quick sidenote, although many of the skeins resulting from a skein computation are not fundamental curves, some can be easily realized as combinations of fundamental curves through simple manipulation, application of isotopy, and Reidemeister moves. Curves like those shown in Figure \ref{fig5} occur so often we give them a variable name as well. This is strictly unnecessary, but very convenient! 

Finally, we note the use of power notation on variables to represent parallel copies of fundamental curves. For example, $x^{2}$ indicates 2 parallel copies of $x$, while $z^{3}$ indicates 3 parallel copies of $z$, etc. Thus, our task of resolving a skein in the genus-2 handlebody amounts to reducing skeins to linear combinations of powers of $x, y, z, y', y''$ and the unknot. 
\noindent Next, we provide a motivational example of how such skein computations are performed for framed skeins in the genus-2 handlebody.

%%%%%%%%%%%%%%%%%%%
\subsection{Example: A Sample Skein Computation} \label{ex1.1}
%%%%%%%%%%%%%%%%%%%
The skein shown in Figure \ref{fig6} can be resolved as a polynomial combination of fundamental curves in $\mathbb{C}[t,t^{-1}]$. In the first row of Figure \ref{fig6} the crossings are resolved one at a time using the skein relations from Figure \ref{fig3}. Resolving these 2 crossings produces $2^2=4$ terms. We note that the first three terms are collections of the fundamental curves given in Figure \ref{fig4} and Figure \ref{fig5}. In the second row of Figure \ref{fig6}, we consider the last term; this skein can be resolved by performing a Reidemeister $R1$ move to induce a crossing by twisting the lower portion of the skein. We present the result at the bottom of Figure \ref{fig6}; we have fully resolved the skein as a linear combination of fundamental curves with coefficients in $\mathbb{C}[t,t^{-1}]$. We can, of course, further simplify the resulting polynomial by distributing and combining like terms; we have produced the skein polynomial associated to the given skein.

\begin{center}
	\begin{figure}[!ht]
		\input{newFig_Ex1a.tikz}
            \input{newFig_Ex1b.tikz}
            \caption{In the first row, we resolve both crossings. The second row induces a crossing by a $R1$ move. The third row expresses the skein as its resulting $\mathbb{C}[t,t^{-1}]$ representation.}
            \label{fig6}
	\end{figure}
\end{center} 

We note here that inducing the crossing in the lower part of the skein was a choice. We could have instead induced the crossing in the upper portion of the skein. In this case, the choice of induction location does not affect our result, but in many cases it does, sometimes negatively.

%%%%%%%%%%%%%%%%%%%%
\subsection{Example: A Problematic Skein Computation} \label{ex1.2}
%%%%%%%%%%%%%%%%%%%%

\begin{center}
	\begin{figure}[!ht]
		\input{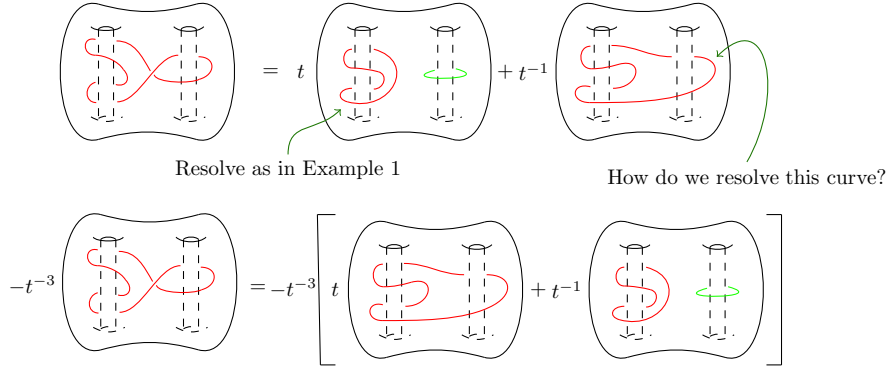}
            \caption{The skein computation described in Example \ref{ex1.2}. A bad choice of crossing location produces the original skein - bad choices compound exponentially! }
            \label{fig7}
	\end{figure}
\end{center}

In this case, the choice of induction location affects the result in a negative manner. Consider the given skein in the upper left of Figure \ref{fig7}. As in Example \ref{ex1.1} we resolve all existing crossings (1 crossing produces $2^1=2$ terms) and determine how to manipulate the skein to produce fundamental curves. The first term is resolved as in Example \ref{ex1.1} (by inducing a crossing on either the upper portion or lower portion by a Reidemeister (R1) move). However, In the second row of Figure \ref{fig7}, we consider the second term; there are multiple places we may induce a crossing. Observe that when we induce a crossing and resolve in the second row of Figure \ref{fig7} we recover the original skein. This was a bad location to induce a crossing! Not only did we end up back at the curve we started with, but we have also introduced another term to the computation.

%%%%%%%%%%%%%%%%%%%%
\subsection{A Method for Skein Computations} \label{sec_method}
%%%%%%%%%%%%%%%%%%%%
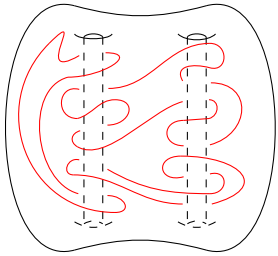
\begin{wrapfigure}{l}[0pt]{5cm}
    \begin{tikzpicture}
	\begin{pgfonlayer}{nodelayer}
		\node [style=none] (79) at (-1.825, 3.5) {};
		\node [style=none] (80) at (1.675, 3.5) {};
		\node [style=none] (81) at (-1.825, -3) {};
		\node [style=none] (82) at (1.675, -3) {};
		\node [style=none] (83) at (-1.825, 2.75) {};
		\node [style=none] (84) at (-0.825, 2.75) {};
		\node [style=none] (85) at (-1.575, 2.675) {};
		\node [style=none] (86) at (-1.075, 2.675) {};
		\node [style=none] (87) at (0.925, 2.75) {};
		\node [style=none] (88) at (1.925, 2.75) {};
		\node [style=none] (89) at (1.175, 2.675) {};
		\node [style=none] (90) at (1.675, 2.675) {};
		\node [style=none] (91) at (-1.825, -2.25) {};
		\node [style=none] (92) at (-0.825, -2.25) {};
		\node [style=none] (93) at (-1.575, -2.3) {};
		\node [style=none] (94) at (-1.075, -2.3) {};
		\node [style=none] (95) at (0.925, -2.25) {};
		\node [style=none] (96) at (1.925, -2.25) {};
		\node [style=none] (97) at (1.175, -2.3) {};
		\node [style=none] (98) at (1.675, -2.3) {};
		\node [style=none] (99) at (-1.7, 1.275) {};
		\node [style=none] (101) at (1.775, 1.85) {};
		\node [style=none] (102) at (1.775, 2.5) {};
		\node [style=none] (103) at (1.075, -0.225) {};
		\node [style=none] (106) at (1.775, -0.2) {};
		\node [style=none] (136) at (1.05, 1.5) {};
		\node [style=none] (137) at (-0.925, 1.3) {};
		\node [style=none] (138) at (1.075, 0.45) {};
		\node [style=none] (140) at (1.8, 1.45) {};
		\node [style=none] (141) at (-1, 2.25) {};
		\node [style=none] (142) at (-2, 1.6) {};
		\node [style=none] (143) at (-1.7, 2.175) {};
		\node [style=none] (144) at (-2.5, 2.5) {};
		\node [style=none] (145) at (-0.75, -2) {};
		\node [style=none] (146) at (-0.95, -1.55) {};
		\node [style=none] (147) at (-1.725, -1.5) {};
		\node [style=none] (148) at (-1.7, -0.6) {};
		\node [style=none] (149) at (-1.85, -1) {};
		\node [style=none] (150) at (1.05, -1.775) {};
		\node [style=none] (151) at (1.825, -1.775) {};
		\node [style=none] (152) at (1, -0.5) {};
		\node [style=none] (153) at (-0.95, -0.85) {};
		\node [style=none] (154) at (1.075, -1.05) {};
		\node [style=none] (155) at (1.775, -0.975) {};
		\node [style=none] (156) at (1.775, -1.5) {};
		\node [style=none] (157) at (-1.75, 0.5) {};
		\node [style=none] (158) at (-1.75, -0.25) {};
		\node [style=none] (159) at (1.75, 0.5) {};
		\node [style=none] (160) at (1.75, 1) {};
		\node [style=none] (161) at (-0.75, 1) {};
		\node [style=none] (162) at (-0.95, 0.325) {};
		\node [style=none] (163) at (1.05, 0.95) {};
		\node [style=none] (164) at (-1.675, 0.175) {};
	\end{pgfonlayer}
	\begin{pgfonlayer}{edgelayer}
		\draw [in=-165, out=-15] (79.center) to (80.center);
		\draw [in=-15, out=15] (80.center) to (82.center);
		\draw [in=15, out=165] (82.center) to (81.center);
		\draw [in=165, out=-165] (81.center) to (79.center);
		\draw [bend right=60, looseness=0.50] (83.center) to (84.center);
		\draw [bend left=45, looseness=0.75] (85.center) to (86.center);
		\draw [bend right=45, looseness=0.75] (87.center) to (88.center);
		\draw [bend left=45, looseness=0.75] (89.center) to (90.center);
		\draw [style=dashed1] (85.center) to (93.center);
		\draw [style=dashed1] (86.center) to (94.center);
		\draw [style=dashed1] (89.center) to (97.center);
		\draw [style=dashed1] (90.center) to (98.center);
		\draw [style=dashed1, bend left] (93.center) to (94.center);
		\draw [style=dashed1, bend left] (97.center) to (98.center);
		\draw [style=dashed1, bend right] (95.center) to (96.center);
		\draw [style=dashed1, bend right] (91.center) to (92.center);
		\draw [style=redline, in=0, out=0, looseness=2.00] (141.center) to (142.center);
		\draw [style=redline, in=45, out=-150, looseness=4.50] (143.center) to (144.center);
		\draw [style=redline, in=-180, out=-135, looseness=1.25] (144.center) to (145.center);
		\draw [style=redline, in=0, out=0, looseness=2.50] (146.center) to (145.center);
		\draw [style=redline, in=-165, out=165] (147.center) to (142.center);
		\draw [style=redline, in=-180, out=165, looseness=2.00] (148.center) to (149.center);
		\draw [style=redline, in=-180, out=0] (149.center) to (150.center);
		\draw [style=redline, in=0, out=15, looseness=2.75] (151.center) to (152.center);
		\draw [style=redline, in=-165, out=-180, looseness=3.25] (152.center) to (154.center);
		\draw [style=redline, in=0, out=15, looseness=2.50] (155.center) to (156.center);
		\draw [style=redline, in=-30, out=-180] (156.center) to (153.center);
		\draw [style=redline, in=0, out=15, looseness=1.75] (106.center) to (140.center);
		\draw [style=redline] (159.center) to (138.center);
		\draw [style=redline, in=15, out=0, looseness=1.50] (160.center) to (159.center);
		\draw [style=redline, in=-180, out=180, looseness=2.25] (103.center) to (138.center);
		\draw [style=redline, in=120, out=165, looseness=1.50] (101.center) to (136.center);
		\draw [style=redline, in=-15, out=0, looseness=2.00] (102.center) to (101.center);
		\draw [style=redline, in=0, out=-180] (102.center) to (137.center);
		\draw [style=redline, in=-180, out=165, looseness=2.00] (99.center) to (157.center);
		\draw [style=redline, in=-180, out=15] (157.center) to (161.center);
		\draw [style=redline, in=0, out=0, looseness=2.25] (161.center) to (162.center);
		\draw [style=redline, in=-30, out=-165, looseness=0.75] (163.center) to (158.center);
		\draw [style=redline, in=150, out=-165, looseness=1.50] (164.center) to (158.center);
	\end{pgfonlayer}
\end{tikzpicture}
    \caption{A more complicated skein: Where do we induce crossings?}
    \label{fig8}
\end{wrapfigure}

Thus far, our method for performing skein computations is to first resolve all the given crossings, cleverly induce new crossings via a series of Reidemeister moves, and then resolve the new crossings we induce. In Example \ref{ex1.1} and \ref{ex1.2}, we see that there is some choice as to where to induce new crossings using Reidemeister moves; some choices are good choices, and result in the resolution of skeins as fundamental (basis) curves. Bad choices result in the recovery of the skein we started with (or an equivalently unknown curve) and the computation produces additional powers and terms which complicate the system even more. Furthermore, given each crossing resolves the skein into 2 terms, bad choices compound exponentially! 

More complicated skeins, like the one shown in Figure \ref{fig8}, introduce even more potential locations to induce a new crossing (after all existing crossings have been resolved). In practice, and in more complicated settings, only expertise and years of experience dictate how a skein should be cleverly isotoped to be realized as a combination of fundamental curves. The process can be frustrating at best for a novice!

In this paper, we propose a series of algorithms to resolve any skein in a genus-2 handlebody as a combination of fundamental curves given in Figure \ref{fig4}, removing the need for clever realization and crossing induction by implementing a sorting method across each strand in the skein projection diagram. This effectively reduces a skein in the genus-2 handlebody $H$ to a skein in a 2-punctured disk $D$. Given that all curves in the $D$ are the desired fundamental curves, this procedure gives us the resolution of any skein in $H$ as basis elements of the Kauffman Bracket Skein Module of $H$. In addition to providing an upper bound on the number of crossings induced (while minimal, we do not guarantee minimum) these algorithms provide a method for automating skein computations (i.e., skeins resolve themselves!) This generalized method removes the need for cleverly realizing more complicated skeins as combinations of known skeins. Additionally, performing all computations with thoughtful reductions at each stage manages the overabundance of terms that increase in computational complexity very efficiently.

In this case, the main difficulty we face in automation is proper input; the current input of our computations are diagrammatic. In order to convert from a diagram to a more algorithmic input (e.g. an array), we must record all the necessary information about our skein system in the given projection. We will establish a notation similar to the format of a braid representation describing the self-intersections of the skein and the relative position of each component to the other components. We then show this representation recovers a unique skein (Lemma \ref{lem1}) and apply Algorithms \ref{alg1}-\ref{alg5} to fully resolve it as a combination of fundamental curves. The result is a complex Laurent polynomial in $x$, $y$, $z$, and the unknot ($t^{-2}+t^{2}$). Note that, although unnecessary, the auxiliary curves $y'$ and $y''$ can be replaced (rather than resolved) by their values in Figure \ref{fig5} for a more efficient computation.

\newpage

%%%%%%%%%%%%%
%%%%%%%%%%%%%%%%%%
\section{Skein Representation: Diagram to Array and Skein Input Parameters}
%%%%%%%%%%%%%%%%%%
%%%%%%%%%%%%%
\subsection{Diagram to Array}
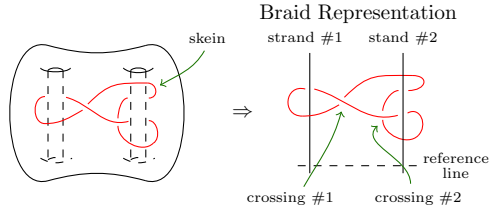
\begin{figure}[!h]
    \begin{tikzpicture}[scale=0.8, {every node/.style}={scale=0.8}]
	\begin{pgfonlayer}{nodelayer}
		\node [style=none] (0) at (-6.55, 2.25) {};
		\node [style=none] (1) at (-3.3, 2.25) {};
		\node [style=none] (2) at (-6.55, -2.25) {};
		\node [style=none] (3) at (-3.3, -2.25) {};
		\node [style=none] (4) at (-6.8, 1.5) {};
		\node [style=none] (5) at (-5.8, 1.5) {};
		\node [style=none] (6) at (-6.55, 1.425) {};
		\node [style=none] (7) at (-6.05, 1.425) {};
		\node [style=none] (8) at (-4.05, 1.5) {};
		\node [style=none] (9) at (-3.05, 1.5) {};
		\node [style=none] (10) at (-3.8, 1.425) {};
		\node [style=none] (11) at (-3.3, 1.425) {};
		\node [style=none] (12) at (-6.8, -1.5) {};
		\node [style=none] (13) at (-5.75, -1.5) {};
		\node [style=none] (14) at (-6.55, -1.55) {};
		\node [style=none] (15) at (-6.05, -1.55) {};
		\node [style=none] (16) at (-4.05, -1.5) {};
		\node [style=none] (17) at (-3, -1.5) {};
		\node [style=none] (18) at (-3.8, -1.55) {};
		\node [style=none] (19) at (-3.3, -1.55) {};
		\node [style=none] (23) at (-5.95, 0.525) {};
		\node [style=none] (25) at (-4.225, -0.425) {};
		\node [style=none] (26) at (-6.65, 0.55) {};
		\node [style=none] (27) at (-4.225, -0.175) {};
		\node [style=none] (29) at (-3.875, -0.275) {};
		\node [style=none] (34) at (-3.05, -1) {};
		\node [style=none] (40) at (-6.65, -0.35) {};
		\node [style=none] (41) at (-3.2, -0.25) {};
		\node [style=none] (140) at (-0.075, 0) {$\Rightarrow$};
		\node [style=none] (367) at (2.15, 2) {};
		\node [style=none] (370) at (2.15, -2) {};
		\node [style=none] (404) at (-5.4, 0.1) {};
		\node [style=none] (405) at (-5.225, 0.3) {};
		\node [style=none] (406) at (-3.25, 1) {};
		\node [style=none] (407) at (-3.9, 0.5) {};
		\node [style=none] (408) at (-3.2, 0.5) {};
		\node [style=none] (409) at (2.3, 0.8) {};
		\node [style=none] (410) at (4.575, -0.225) {};
		\node [style=none] (411) at (2, 0.8) {};
		\node [style=none] (412) at (4.575, 0.075) {};
		\node [style=none] (413) at (5.075, -0.075) {};
		\node [style=none] (414) at (5.75, -0.75) {};
		\node [style=none] (415) at (1.9, -0.1) {};
		\node [style=none] (416) at (5.4, 0) {};
		\node [style=none] (417) at (3.125, 0.275) {};
		\node [style=none] (418) at (3.3, 0.525) {};
		\node [style=none] (419) at (5.55, 1.25) {};
		\node [style=none] (420) at (5.1, 0.75) {};
		\node [style=none] (421) at (5.425, 0.75) {};
		\node [style=none] (422) at (5.25, 2) {};
		\node [style=none] (423) at (5.25, -2) {};
		\node [style=none] (424) at (-2.825, 0.925) {};
		\node [style=none] (425) at (-1.35, 2.075) {};
		\node [style=none] (426) at (-1.25, 2.5) {};
		\node [style=none] (427) at (-1.25, 2.5) {\scriptsize{skein}};
		\node [style=none] (428) at (1.75, -1.75) {};
		\node [style=none] (429) at (5.75, -1.75) {};
		\node [style=none] (430) at (2, 2.5) {\scriptsize{strand \#1}};
		\node [style=none] (431) at (5.25, 2.5) {\scriptsize{stand \#2}};
		\node [style=none] (432) at (3.75, 3.3) {Braid Representation};
		\node [style=none] (433) at (3.2, 0.05) {};
		\node [style=none] (434) at (1.95, -2.45) {};
		\node [style=none] (435) at (4.4, -0.425) {};
		\node [style=none] (436) at (5.5, -2.35) {};
		\node [style=none] (437) at (1.5, -2.85) {\scriptsize{crossing \#1}};
		\node [style=none] (438) at (5.75, -2.85) {\scriptsize{crossing \#2}};
		\node [style=none] (439) at (7, -1.5) {\scriptsize{reference}};
		\node [style=none] (440) at (7, -2) {\scriptsize{line}};
	\end{pgfonlayer}
	\begin{pgfonlayer}{edgelayer}
		\draw [in=-165, out=-15] (0.center) to (1.center);
		\draw [in=-15, out=15] (1.center) to (3.center);
		\draw [in=15, out=165] (3.center) to (2.center);
		\draw [in=165, out=-165] (2.center) to (0.center);
		\draw [bend right=60, looseness=0.50] (4.center) to (5.center);
		\draw [bend left=45, looseness=0.75] (6.center) to (7.center);
		\draw [bend right=45, looseness=0.75] (8.center) to (9.center);
		\draw [bend left=45, looseness=0.75] (10.center) to (11.center);
		\draw [style=dashed1] (6.center) to (14.center);
		\draw [style=dashed1] (7.center) to (15.center);
		\draw [style=dashed1] (10.center) to (18.center);
		\draw [style=dashed1] (11.center) to (19.center);
		\draw [style=dashed1, bend left] (14.center) to (15.center);
		\draw [style=dashed1, bend left] (18.center) to (19.center);
		\draw [style=dashed1, bend right] (16.center) to (17.center);
		\draw [style=dashed1, bend right] (12.center) to (13.center);
		\draw (367.center) to (370.center);
		\draw [style=redline, in=-165, out=-15] (23.center) to (29.center);
		\draw [style=redline, in=30, out=0] (41.center) to (34.center);
		\draw [style=redline, in=-135, out=-90] (25.center) to (34.center);
		\draw [style=redline, in=180, out=-180, looseness=1.25] (26.center) to (40.center);
		\draw [style=redline, in=-150, out=0] (40.center) to (404.center);
		\draw [style=redline, in=180, out=45, looseness=1.25] (405.center) to (406.center);
		\draw [style=redline, in=-180, out=90, looseness=0.75] (27.center) to (407.center);
		\draw [style=redline, in=0, out=0, looseness=1.75] (406.center) to (408.center);
		\draw [style=redline, in=-165, out=-15] (409.center) to (413.center);
		\draw [style=redline, in=45, out=30, looseness=1.50] (416.center) to (414.center);
		\draw [style=redline, in=-135, out=-90] (410.center) to (414.center);
		\draw [style=redline, in=180, out=165, looseness=2.00] (411.center) to (415.center);
		\draw [style=redline, in=-150, out=0] (415.center) to (417.center);
		\draw [style=redline, in=-180, out=45, looseness=1.25] (418.center) to (419.center);
		\draw [style=redline, in=-180, out=90, looseness=0.75] (412.center) to (420.center);
		\draw [style=redline, in=0, out=0, looseness=3.00] (419.center) to (421.center);
		\draw (422.center) to (423.center);
		\draw [style=arrow4, in=-120, out=30] (424.center) to (425.center);
		\draw [style=dashed1] (428.center) to (429.center);
		\draw [style=arrow4, in=90, out=-135] (435.center) to (436.center);
		\draw [style=arrow4, in=60, out=-105] (433.center) to (434.center);
	\end{pgfonlayer}
\end{tikzpicture}
    \caption{Braid Representation}
    \label{fig9}
\end{figure}

 We start by converting our skein system to a braid diagram as shown in Figure \ref{fig9}. Note that nothing about the skein itself has changed, and we view the holes in the genus-2 handlebody as fixed strands of an underlying base braid. When we perform skein computations, the strands of the base braid stay fixed; crossings between the skein and these strands need not be resolved - only crossings between the skein and itself (self-intersections) are resolved using the given skein relations.
 
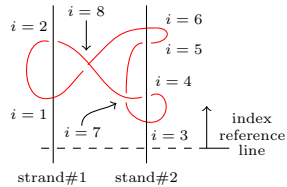
\begin{figure}[!h]
    \begin{tikzpicture}[scale=0.8, {every node/.style}={scale=0.8}]
	\begin{pgfonlayer}{nodelayer}
		\node [style=none] (367) at (-1.6, 2.5) {};
		\node [style=none] (370) at (-1.6, -2.75) {};
		\node [style=none] (409) at (-1.45, 1.3) {};
		\node [style=none] (410) at (0.825, -0.725) {};
		\node [style=none] (411) at (-1.75, 1.3) {};
		\node [style=none] (412) at (0.825, -0.425) {};
		\node [style=none] (413) at (1.35, -0.55) {};
		\node [style=none] (414) at (2, -1.25) {};
		\node [style=none] (415) at (-1.85, -0.6) {};
		\node [style=none] (416) at (1.65, -0.5) {};
		\node [style=none] (417) at (-0.6, 0.325) {};
		\node [style=none] (418) at (-0.45, 0.525) {};
		\node [style=none] (419) at (1.8, 1.75) {};
		\node [style=none] (420) at (1.375, 1.25) {};
		\node [style=none] (421) at (1.65, 1.25) {};
		\node [style=none] (422) at (1.5, 2.5) {};
		\node [style=none] (423) at (1.5, -2.75) {};
		\node [style=none] (428) at (-2, -2.25) {};
		\node [style=none] (429) at (3, -2.25) {};
		\node [style=none] (430) at (-1.625, -3.25) {\scriptsize{strand\#1}};
		\node [style=none] (431) at (1.5, -3.25) {\scriptsize{stand\#2}};
		\node [style=none] (433) at (-2.4, -1.125) {\scriptsize{$i=1$}};
		\node [style=none] (434) at (-2.4, 1.8) {\scriptsize{$i=2$}};
		\node [style=none] (435) at (0.65, -0.925) {};
		\node [style=none] (436) at (3.5, -0.85) {};
		\node [style=none] (439) at (5, -1.8) {\scriptsize{reference}};
		\node [style=none] (440) at (5, -2.3) {\scriptsize{line}};
		\node [style=none] (441) at (3.25, -2.25) {};
		\node [style=none] (442) at (4.25, -2.25) {};
		\node [style=none] (443) at (3.5, -2.25) {};
		\node [style=none] (444) at (5, -1.25) {};
		\node [style=none] (445) at (5, -1.25) {\scriptsize{index}};
		\node [style=none] (446) at (2.75, 2.05) {\scriptsize{$i=6$}};
		\node [style=none] (447) at (-0.5, 2.3) {\scriptsize{$i=8$}};
		\node [style=none] (448) at (-0.5, 2) {};
		\node [style=none] (449) at (-0.5, 1) {};
		\node [style=none] (450) at (2.75, 1.05) {\scriptsize{$i=5$}};
		\node [style=none] (451) at (2.4, -0.025) {\scriptsize{$i=4$}};
		\node [style=none] (452) at (2.275, -1.8) {\scriptsize{$i=3$}};
		\node [style=none] (453) at (-0.625, -1.725) {\scriptsize{$i=7$}};
		\node [style=none] (454) at (-0.65, -1.5) {};
		\node [style=none] (455) at (0.5, -0.75) {};
	\end{pgfonlayer}
	\begin{pgfonlayer}{edgelayer}
		\draw (367.center) to (370.center);
		\draw [style=redline, in=-165, out=-15] (409.center) to (413.center);
		\draw [style=redline, in=45, out=30, looseness=1.50] (416.center) to (414.center);
		\draw [style=redline, in=-135, out=-90] (410.center) to (414.center);
		\draw [style=redline, in=-180, out=165, looseness=1.25] (411.center) to (415.center);
		\draw [style=redline, in=-135, out=0] (415.center) to (417.center);
		\draw [style=redline, in=-180, out=45, looseness=1.25] (418.center) to (419.center);
		\draw [style=redline, in=-180, out=90, looseness=0.75] (412.center) to (420.center);
		\draw [style=redline, in=0, out=0, looseness=3.00] (419.center) to (421.center);
		\draw (422.center) to (423.center);
		\draw [style=dashed1] (428.center) to (429.center);
		\draw (441.center) to (442.center);
		\draw [style=arrow3] (436.center) to (443.center);
		\draw [style=arrow3] (449.center) to (448.center);
		\draw [style=arrow3, in=90, out=-150] (455.center) to (454.center);
	\end{pgfonlayer}
\end{tikzpicture}
    \caption{Indexing a single component}
    \label{fig10}
\end{figure}

\begin{figure}[!h]
    \begin{tikzpicture}[scale=0.8, {every node/.style}={scale=0.8}]
	\begin{pgfonlayer}{nodelayer}
		\node [style=none] (367) at (-1.85, 3.75) {};
		\node [style=none] (370) at (-1.85, -3.25) {};
		\node [style=none] (409) at (-1.7, 1.05) {};
		\node [style=none] (411) at (-2, 1.05) {};
		\node [style=none] (412) at (1.625, -0.675) {};
		\node [style=none] (413) at (2.1, -1.8) {};
		\node [style=none] (414) at (2.1, -1.8) {};
		\node [style=none] (415) at (-2.1, -0.6) {};
		\node [style=none] (416) at (1.9, -0.75) {};
		\node [style=none] (417) at (-0.275, 0.125) {};
		\node [style=none] (418) at (-0.125, 0.275) {};
		\node [style=none] (419) at (1.55, 1.25) {};
		\node [style=none] (420) at (1.875, 1.25) {};
		\node [style=none] (422) at (1.75, 3.75) {};
		\node [style=none] (423) at (1.75, -3.25) {};
		\node [style=none] (430) at (-1.875, -3.75) {\scriptsize{strand\#1}};
		\node [style=none] (431) at (1.75, -3.75) {\scriptsize{stand\#2}};
		\node [style=none] (433) at (-1.025, -3.1) {\scriptsize{$i=1$}};
		\node [style=none] (434) at (-0.85, -2.45) {\scriptsize{$i=2$}};
		\node [style=none] (435) at (0.625, -0.625) {};
		\node [style=none] (436) at (-0.075, -0.275) {};
		\node [style=none] (439) at (5, -3.05) {\scriptsize{crossing\#2}};
		\node [style=none] (444) at (3, -1.25) {};
		\node [style=none] (445) at (2.5, 3.35) {\scriptsize{crossing\#1}};
		\node [style=none] (446) at (-1.125, 1.5) {\scriptsize{$i=6$}};
		\node [style=none] (447) at (-1.125, 2.725) {\scriptsize{$i=8$}};
		\node [style=none] (448) at (-2.25, 2.25) {};
		\node [style=none] (449) at (-1.7, 2) {};
		\node [style=none] (450) at (-1.2, 0.05) {\tiny{$i=5$}};
		\node [style=none] (451) at (-1.2, -1.15) {\tiny{$i=4$}};
		\node [style=none] (452) at (-0.9, -1.875) {\scriptsize{$i=3$}};
		\node [style=none] (453) at (-0.875, 2.125) {\scriptsize{$i=7$}};
		\node [style=none] (456) at (-2, -2.75) {};
		\node [style=none] (457) at (-1.95, 1.85) {};
		\node [style=none] (458) at (-2.25, -1) {};
		\node [style=none] (459) at (0.25, 0) {};
		\node [style=none] (460) at (0.475, -0.575) {};
		\node [style=none] (461) at (0.325, -0.75) {};
		\node [style=none] (462) at (-1.75, -1.75) {};
		\node [style=none] (463) at (-2.05, -1.85) {};
		\node [style=none] (464) at (-1.675, -2.675) {};
		\node [style=none] (465) at (2.6, -2.175) {\scriptsize{$i=10$}};
		\node [style=none] (466) at (2.65, -0.45) {\scriptsize{$i=11$}};
		\node [style=none] (467) at (2.625, 0.35) {\scriptsize{$i=12$}};
		\node [style=none] (468) at (2.75, 1.6) {\scriptsize{$i=13$}};
		\node [style=none] (469) at (0.25, 3.35) {\scriptsize{$i=14$}};
		\node [style=none] (470) at (0.25, 3) {};
		\node [style=none] (471) at (-0.1, 0.5) {};
		\node [style=none] (472) at (2.25, -3) {};
		\node [style=none] (473) at (2.9, -3.05) {\scriptsize{$i=9$}};
		\node [style=none] (474) at (5.775, -1.5) {\scriptsize{crossing\#3}};
		\node [style=none] (475) at (3.45, -1.5) {\scriptsize{$i=15$}};
	\end{pgfonlayer}
	\begin{pgfonlayer}{edgelayer}
		\draw (367.center) to (370.center);
		\draw [style=redline, in=180, out=0] (409.center) to (413.center);
		\draw [style=redline, in=15, out=-15, looseness=1.75] (416.center) to (414.center);
		\draw [style=redline, in=-180, out=-180, looseness=1.50] (411.center) to (415.center);
		\draw [style=redline, in=-135, out=0] (415.center) to (417.center);
		\draw [style=redline, in=-165, out=45, looseness=1.25] (418.center) to (419.center);
		\draw [style=redline, in=0, out=150, looseness=2.25] (412.center) to (420.center);
		\draw (422.center) to (423.center);
		\draw [style=redline, in=60, out=30, looseness=1.75] (448.center) to (449.center);
		\draw [style=redline, in=180, out=-150] (448.center) to (456.center);
		\draw [style=redline, in=-180, out=-150, looseness=1.25] (457.center) to (458.center);
		\draw [style=redline, in=-165, out=0, looseness=0.50] (458.center) to (459.center);
		\draw [style=redline, in=15, out=45, looseness=2.25] (460.center) to (459.center);
		\draw [style=redline, in=-135, out=15, looseness=0.75] (462.center) to (461.center);
		\draw [style=redline, in=45, out=-150, looseness=2.00] (463.center) to (464.center);
		\draw [style=arrow3, in=-75, out=60, looseness=1.25] (471.center) to (470.center);
		\draw [style=arrow3, in=135, out=-135, looseness=0.75] (436.center) to (472.center);
		\draw [style=arrow3] (435.center) to (444.center);
	\end{pgfonlayer}
\end{tikzpicture}
    \caption{Indexing multiple components}
    \label{fig11}
\end{figure}
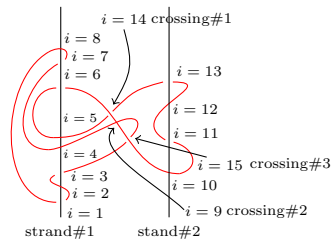

 Next, we define a reference line at the bottom of the braid diagram (Figure \ref{fig10}) and index the crossings between the skein and the base braid from bottom to top along each strand of the braid. We will also index crossings between the skein and itself. The indexing sets between strands 1, 2, and the set of skein self-intersections do not have to be related, and, using a different indexing scheme along each strand will not affect the final result; it is only important that elements along each strand are ordered from bottom to top. Skeins with multiple components, Figure \ref{fig11}, are indexed along each strand from bottom to top regardless of to which component that portion of the skein belongs. We will identify each of the crossings in the skein by labeling them on the diagram. The order of the crossings determines the order in which they will be resolved. Note that this ordering is inherently different than the index ordering we assigned. To accentuate this difference, we will label the crossings with rational numbers $(0.01, 0.02, ..., 0.0n)$, not to be confused with the counting numbers we use for indices based on our reference line.

In general, we can label our crossings by $1 \times 10^{-k_{r}}, 2 \times 10^{-k_{r}}, ..., r \times 10^{-k_{r}}$ where $k_{r}$ is the magnitude of the number of crossings. For example, a skein with 1-99 crossings are labeled as 0.01, 0.02, ..., 0.99, a skein with 100-999 crossings are labeled 0.001, 0.002,...,0.999, etc. Finally, at the locations where the skein intersects the strands of the base braid we will note whether the skein goes over ($+n$) or under ($-n$) the base braid $n$. With this information and the labeling on the set of skein self-intersections we can form an array, $E^{l}$, for each component $l$ of the skein. Among all components of the skein if there are $s$ self-intersections and $r$ intersections between the skein and the braid then there are $r+2s$ points to record. i.e., the label on each self-intersection is recorded twice, once for the overcrossing and once for the undercrossing, although these over/undercrossings may be in the same or different components of the skein.

We outline how to record skein information at every point along the skein. The points we consider are locations where the skein crosses the strands of the base braid in the given projection. In each component $l$ of the skein $K$, and at each point $p$, we note the skein intersection information ($e^l_p$), skein index information ($i^l_p$), and skein orientation information ($q^l_p$). Recording this information over every component, we prove (Lemma \ref{lem1}) that this representation $(E,I,Q)$ uniquely recovers the skein $K$

%%%%%%%%%%%%%%%
\subsection{Crossing Input Information (E)}
%%%%%%%%%%%%%%%

In each component $l$ of the skein $K$, let $p_{l}$ be the number of points of interest (places where the skein intersects the strands of the base braid or skein self-intersections). If there are $k$ total components in the skein, for each component $0 < l \leq k$, we form an array:

\begin{equation*}
    E^{l}=[e_{0}^{l},e_{1}^{l}, ..., e_{p_{l}-1}^{l}]
\end{equation*}

where $e_{j}^{l}$, $0 \leq j \leq (p_{l}-1)$ is determined by the following, letting $S$ be the set of skein self-intersections:

\begin{table}[h]
\begin{center}
\begin{tabular}{ |c|c|}
\hline
Skein Intersection Input Information & Value:  \\
\hline
If the skein crosses over strand 1 of the base braid: & $e_{j}^{l}=1$  \\ 
If the skein crosses under strand 1 of the base braid: & $e_{j}^{l}=-1$  \\ 
If the skein crosses over strand 2 of the base braid: & $e_{j}^{l}=2$  \\ 
If the skein crosses under strand 2 of the base braid: & $e_{j}^{l}=-2$  \\ 
\hline
If the skein crosses over a self-intersection ($s \in S$) & $e_{j}^{l}=s$  \\ 
If the skein crosses under a self-intersection ($s \in S$) & $e_{j}^{l}=-s$  \\ 
\hline
\end{tabular}
\end{center}
\end{table}
\vspace{0.2cm}
We record $E^{l}$ for each component $l$, $0 < l \leq k$, and from this we form the skein intersection array:

\begin{equation*}
    E=[E^{1},E^{2},...,E^{l},...,E^{k}]
\end{equation*}

\noindent where $k$ is the number of components in the skein $K$.

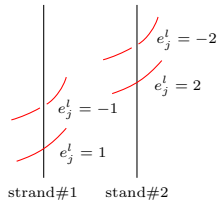
\begin{figure}[!h]
    \centering
    \begin{tikzpicture}[scale=0.8, {every node/.style}={scale=0.8}]
	\begin{pgfonlayer}{nodelayer}
		\node [style=none] (367) at (-1.6, 2.5) {};
		\node [style=none] (370) at (-1.6, -3.25) {};
		\node [style=none] (411) at (-2.5, -2.7) {};
		\node [style=none] (415) at (-0.85, -1.6) {};
		\node [style=none] (420) at (-0.85, 0.1) {};
		\node [style=none] (422) at (1.5, 2.5) {};
		\node [style=none] (423) at (1.5, -3.25) {};
		\node [style=none] (430) at (-1.625, -3.75) {\scriptsize{strand\#1}};
		\node [style=none] (431) at (1.5, -3.75) {\scriptsize{stand\#2}};
		\node [style=none] (456) at (-2.675, -1.325) {};
		\node [style=none] (457) at (-1.725, -0.925) {};
		\node [style=none] (459) at (0.675, -0.45) {};
		\node [style=none] (460) at (2.325, 0.65) {};
		\node [style=none] (464) at (-1.5, -0.8) {};
		\node [style=none] (468) at (-0.3, -2.45) {\scriptsize{$e^{l}_{j}=1$}};
		\node [style=none] (469) at (2.85, -0.2) {\scriptsize{$e^{l}_{j}=2$}};
		\node [style=none] (471) at (-0.125, -0.975) {\scriptsize{$e^{l}_{j}=-1$}};
		\node [style=none] (473) at (2.2, 2.1) {};
		\node [style=none] (474) at (0.375, 0.675) {};
		\node [style=none] (475) at (1.325, 1.075) {};
		\node [style=none] (476) at (1.55, 1.2) {};
		\node [style=none] (477) at (3.15, 1.35) {\scriptsize{$e^{l}_{j}=-2$}};
	\end{pgfonlayer}
	\begin{pgfonlayer}{edgelayer}
		\draw (367.center) to (370.center);
		\draw [style=redline, bend right=15] (411.center) to (415.center);
		\draw (422.center) to (423.center);
		\draw [style=redline, in=-150, out=15, looseness=0.75] (456.center) to (457.center);
		\draw [style=redline, bend right=15] (459.center) to (460.center);
		\draw [style=redline, in=-105, out=30, looseness=0.75] (464.center) to (420.center);
		\draw [style=redline, in=-150, out=15, looseness=0.75] (474.center) to (475.center);
		\draw [style=redline, in=-105, out=30, looseness=0.75] (476.center) to (473.center);
	\end{pgfonlayer}
\end{tikzpicture}
    \caption{The value of entries in the skein intersection array}
    \label{fig12}
\end{figure}

We note that we must choose a reference point and direction to list elements in each component of the skein. That is, we choose an element to begin our array arbitrarily. This choice does not affect our final answer. Figure \ref{fig12} shows how to label skein intersection information every time the skein crosses over or under the base braid (that is, passes in front or behind the hole in the genus-2 surface from our given projection).

%%%%%%%%%%%%%%%
\subsection{Index Input Information (I)}
%%%%%%%%%%%%%%%

Next, we form the index array for each component $l$, starting with the point of interest corresponding to $e_{0}^{l}$, (the first entry in $E^{l}$) after our assigned reference point. We proceed in the same order as the elements of $E^{l}$ so that the entries of $I^{l}$ correspond to the entries of $E^{l}$. We record $i^{l}_{j}$ point by point as assigned relative to our line of reference. For each component $0 < l \leq k$: 

\begin{equation*}
    I^{l}=[i_{0}^{l},i_{1}^{l}, ..., i_{p_{l}-1}^{l}]
\end{equation*}

\begin{equation*}
    I=[I^{1},I^{2},...,I^{l},...,I^{k}]
\end{equation*}

where $k$ is the number of components in the skein $K$.

%%%%%%%%%%%%%%
\subsection{Orientation Input Information (Q)}
%%%%%%%%%%%%%%

We form an array for the orientation of each component:

\begin{equation*}
    Q^{l}=[q_{0}^{l},q_{1}^{l}, ..., q_{p_{l}-1}^{l}]
\end{equation*}

While it is only necessary to define whether the component is traversed clockwise or counterclockwise from the assigned reference point (i.e., the assignment of a single value for each component), we have to make a local decision about orientation when comparing any two elements. Thus, we proceed to assign each element $e_{j}^{l},i_{j}^{l}$ a value $q_{j}^{l}$, $0 \leq j \leq (p_{l}-1)$, corresponding to the location of that point in the given projection, and compare these values between two points to determine a motion from left to right or right to left.

\vspace{0.2cm}
\begin{table}[h]
\begin{center}
\begin{tabular}{ |c|c|} 
\hline
Skein Orientation Input Information & Value:  \\
\hline
If $e_{j}^{l}=1$, and the skein goes from left to right: & $q_{j}^{l}=3$  \\ 
If $e_{j}^{l}=1$, and the skein goes from right to left: & $q_{j}^{l}=4$  \\ 
If $e_{j}^{l}=2$, and the skein goes from left to right: & $q_{j}^{l}=4$  \\ 
If $e_{j}^{l}=2$, and the skein goes from left to right: & $q_{j}^{l}=5$  \\ 
\hline
If $e_{j}^{l}=\pm0.0k$: & $q_{j}^{l}=0$  \\ 
\hline
\end{tabular}
\end{center}
\end{table}
\vspace{0.2cm}

\noindent We record $Q^{l}$ for each component $l$, $0 < l \leq k$, and from this we form the skein intersection array:

\begin{equation*}
    Q=[Q^{1},Q^{2},,...Q^{l},...,Q^{k}]
\end{equation*}

where $k$ is the number of components in the skein $K$.
\subsection{Self-intersections (U)} \label{u}

Finally, when performing skein computations for oriented skeins, we record the sign on each crossing. To determine the sign of the crossing, we look at the orientation of the overcrossing and then undercrossing. We keep track of the signs on the crossings to be resolved because the resolution depends on whether the crossing is positive or negative (as defined by the right hand rule for skeins). We record the signs of the crossings at each stage in an array $U$. When a crossing is resolved, the corresponding entry of $U$ is deleted. 

When we resolve a crossing, we smooth it in both directions (as shown in Figure \ref{fig3}); this smoothing reverses the orientation of one of the sections of the skein. The sign of any other self-intersection (yet to be resolved) that lies along that section may be affected. If both components (the overcrossing and undercrossing) of another crossing lie on this segment, the orientation of both changes, and no sign change is necessary. However, if only one of the components (either the overcrossing or the undercrossing) of another crossing lies on this segment, the sign of that crossing will change in the next stage. Thus, for each new skein term, we record the signs of the remaining crossings for that term with the other skein information.
 
Thus, we store the signs of crossings to be resolved in an array, called "U", as follows:
\begin{equation*}
    U=[u_1,u_2,...u_r]
\end{equation*}

where $r$ is the number of self-intersections in the original skein. As we resolve these crossings, we delete their corresponding entries from $U$; when the skein $K$ is fully resolved, $U=\emptyset$.

\hspace{0.5cm}

We pause the narrative for two examples demonstrating how to record this information for a skein with one component \ref{ex_single} and a skein with two components \ref{ex_multi}. We then summerize the skein input array in \ref{skein_input}.
\newpage

%%%%%%%%%%%%%%
\subsection{Example (Single Component): Crossing, Index, and Orientation Inputs} \label{ex_single}
%%%%%%%%%%%%%%

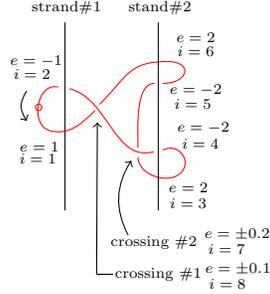
\begin{figure}[!h]
    \centering
        \begin{tikzpicture}[scale=0.8, {every node/.style}={scale=0.8}]
	\begin{pgfonlayer}{nodelayer}
		\node [style=none] (367) at (-1.6, 3.25) {};
		\node [style=none] (370) at (-1.6, -3) {};
		\node [style=none] (409) at (-1.45, 1.05) {};
		\node [style=none] (410) at (0.825, -1.225) {};
		\node [style=none] (411) at (-1.825, 1.1) {};
		\node [style=none] (412) at (0.825, -0.925) {};
		\node [style=none] (413) at (1.35, -1.05) {};
		\node [style=none] (414) at (2.25, -1.75) {};
		\node [style=none] (415) at (-2.1, -0.35) {};
		\node [style=none] (416) at (1.65, -1) {};
		\node [style=none] (417) at (-0.6, 0.325) {};
		\node [style=none] (418) at (-0.45, 0.525) {};
		\node [style=none] (419) at (1.8, 1.95) {};
		\node [style=none] (420) at (1.375, 1.225) {};
		\node [style=none] (421) at (1.65, 1.225) {};
		\node [style=none] (422) at (1.5, 3.25) {};
		\node [style=none] (423) at (1.5, -3) {};
		\node [style=none] (430) at (-1.55, 3.75) {\scriptsize{strand\#1}};
		\node [style=none] (431) at (1.55, 3.75) {\scriptsize{stand\#2}};
		\node [style=none] (433) at (-2.5, -1.275) {\scriptsize{$i=1$}};
		\node [style=none] (434) at (-2.7, 1.525) {\scriptsize{$i=2$}};
		\node [style=none] (446) at (2.75, 2.3) {\scriptsize{$i=6$}};
		\node [style=none] (450) at (2.65, 0.55) {\scriptsize{$i=5$}};
		\node [style=none] (451) at (2.9, -0.775) {\scriptsize{$i=4$}};
		\node [style=none] (457) at (-2.875, 0.975) {};
		\node [style=none] (458) at (-2.875, -0.025) {};
		\node [style=none] (459) at (-2.55, 1.975) {\scriptsize{$e=-1$}};
		\node [style=none] (460) at (-2.45, -0.9) {\scriptsize{$e=1$}};
		\node [style=none] (461) at (2.5, -2.775) {\scriptsize{$i=3$}};
		\node [style=none] (462) at (2.75, 2.75) {\scriptsize{$e=2$}};
		\node [style=none] (463) at (2.75, 1) {\scriptsize{$e=-2$}};
		\node [style=none] (464) at (3, -0.25) {\scriptsize{$e=-2$}};
		\node [style=none] (465) at (2.5, -2.25) {\scriptsize{$e=2$}};
		\node [style=none] (466) at (1.35, -4.075) {\scriptsize{crossing \#2}};
		\node [style=none] (467) at (4.125, -3.8) {\scriptsize{$e=\pm0.2$}};
		\node [style=none] (468) at (3.775, -4.3) {\scriptsize{$i=7$}};
		\node [style=none] (469) at (1.5, -5.125) {\scriptsize{crossing \#1}};
		\node [style=none] (470) at (4.15, -4.95) {\scriptsize{$e=\pm0.1$}};
		\node [style=none] (471) at (3.8, -5.45) {\scriptsize{$i=8$}};
		\node [style=none] (472) at (0.6, -1.35) {};
		\node [style=none] (473) at (0.5, -3.825) {};
		\node [style=none] (474) at (-0.525, -0.075) {};
		\node [style=none] (475) at (-0.55, -5.125) {};
		\node [style=none] (476) at (0, -5.125) {};
		\node [style=reddot1] (477) at (-2.475, 0.425) {};
	\end{pgfonlayer}
	\begin{pgfonlayer}{edgelayer}
		\draw (367.center) to (370.center);
		\draw [style=redline, in=-165, out=-15] (409.center) to (413.center);
		\draw [style=redline, in=45, out=30, looseness=1.50] (416.center) to (414.center);
		\draw [style=redline, in=-135, out=-90] (410.center) to (414.center);
		\draw [style=redline, in=165, out=165, looseness=1.25] (411.center) to (415.center);
		\draw [style=redline, in=-135, out=-15] (415.center) to (417.center);
		\draw [style=redline, in=-180, out=45, looseness=1.25] (418.center) to (419.center);
		\draw [style=redline, in=-180, out=90, looseness=0.75] (412.center) to (420.center);
		\draw [style=redline, in=0, out=0, looseness=3.00] (419.center) to (421.center);
		\draw (422.center) to (423.center);
		\draw [style=arrow2, bend right, looseness=1.25] (457.center) to (458.center);
		\draw [style=arrow2, bend left] (473.center) to (472.center);
		\draw [style=arrow2] (475.center) to (474.center);
		\draw (476.center) to (475.center);
	\end{pgfonlayer}
\end{tikzpicture}
        \caption{Single Component Crossing, Index, and Orientation Inputs}
        \label{fig13}
\end{figure}

The skein in Figure \ref{fig13} has only one component ($l=1$). We set a reference point, and record the skein intersection information $E^{1}=[1,-0.1,2,-2,-0.2,2,-2,0.2,0.1,-1]$ and $E=[E^1]$. From the reference point shown, the index information is recorded $I^1=[1,8,6,5,7,3,4,7,8,2]$, and $I=[I^1]$. For orientation, we record $Q^1=[3,0,4,5,0,4,5,0,0,4]$, and $Q=[Q^1]$. (Note that $dim(E^l)=dim(I^l)=dim(Q^l)$, since $e_{j}^{l}$ corresponds one-to-one with $i_{j}^{l}$ and $q_{j}^{l}$, and $dim(E)=dim(I)=dim(Q)$.) There are 2 crossings to resolve, labeled 0.1 and 0.2. The array $U$ stores the signs of these crossings: $U=[-1,1]$.

%%%%%%%%%%%%%%%
\subsection{Example (Multi-Component): Crossing, Index, and Orientation Inputs} \label{ex_multi}
%%%%%%%%%%%%%%%

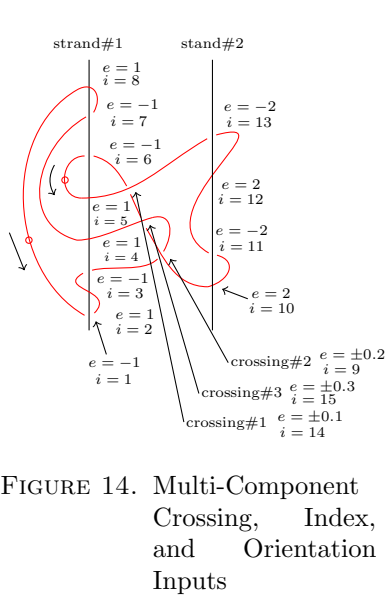
\begin{wrapfigure}{L}[0pt]{5cm}
     \centering
        \begin{tikzpicture}[scale=0.8, {every node/.style}={scale=0.8}]
	\begin{pgfonlayer}{nodelayer}
		\node [style=none] (367) at (-2.1, 4.75) {};
		\node [style=none] (370) at (-2.1, -4.25) {};
		\node [style=none] (409) at (-1.95, 1.55) {};
		\node [style=none] (411) at (-2.25, 1.55) {};
		\node [style=none] (412) at (1.875, -1.675) {};
		\node [style=none] (413) at (-0.725, 0.275) {};
		\node [style=none] (414) at (1.85, -2.8) {};
		\node [style=none] (415) at (-2.35, 0.15) {};
		\node [style=none] (416) at (2.15, -1.75) {};
		\node [style=none] (417) at (1.825, 2.125) {};
		\node [style=none] (418) at (-0.375, -0.5) {};
		\node [style=none] (419) at (1.55, -6.35) {};
		\node [style=none] (420) at (2.125, 2.25) {};
		\node [style=none] (422) at (2, 4.75) {};
		\node [style=none] (423) at (2, -4.25) {};
		\node [style=none] (430) at (-2.125, 5.25) {\scriptsize{strand\#1}};
		\node [style=none] (431) at (2, 5.25) {\scriptsize{stand\#2}};
		\node [style=none] (433) at (-1.275, -5.875) {\scriptsize{$i=1$}};
		\node [style=none] (434) at (-0.6, -4.2) {\scriptsize{$i=2$}};
		\node [style=none] (435) at (2.525, -5.325) {};
		\node [style=none] (436) at (-0.2, -0.85) {};
		\node [style=none] (439) at (3.95, -5.3) {\scriptsize{crossing\#2}};
		\node [style=none] (445) at (2.475, -7.35) {\scriptsize{crossing\#1}};
		\node [style=none] (446) at (-0.625, 1.45) {\scriptsize{$i=6$}};
		\node [style=none] (447) at (-1.025, 4.075) {\scriptsize{$i=8$}};
		\node [style=none] (448) at (-2.5, 3.75) {};
		\node [style=none] (449) at (-1.95, 3) {};
		\node [style=none] (450) at (-1.375, -0.6) {\tiny{$i=5$}};
		\node [style=none] (451) at (-1.025, -1.8) {\tiny{$i=4$}};
		\node [style=none] (452) at (-0.9, -3.025) {\scriptsize{$i=3$}};
		\node [style=none] (453) at (-0.775, 2.725) {\scriptsize{$i=7$}};
		\node [style=none] (456) at (-2.25, -3.75) {};
		\node [style=none] (457) at (-2.2, 2.85) {};
		\node [style=none] (458) at (-2.5, -1.25) {};
		\node [style=none] (459) at (0, -0.5) {};
		\node [style=none] (460) at (0.375, -1.625) {};
		\node [style=none] (461) at (0.175, -1.875) {};
		\node [style=none] (462) at (-2, -2.25) {};
		\node [style=none] (463) at (-2.3, -2.35) {};
		\node [style=none] (464) at (-1.925, -3.675) {};
		\node [style=none] (465) at (3.975, -3.5) {\scriptsize{$i=10$}};
		\node [style=none] (466) at (2.975, -1.45) {\scriptsize{$i=11$}};
		\node [style=none] (467) at (2.95, 0.1) {\scriptsize{$i=12$}};
		\node [style=none] (468) at (3.225, 2.675) {\scriptsize{$i=13$}};
		\node [style=none] (469) at (5, -7.65) {\scriptsize{$i=14$}};
		\node [style=none] (470) at (-0.075, -0.75) {};
		\node [style=none] (471) at (-0.85, 0.55) {};
		\node [style=none] (472) at (0.6, -1.875) {};
		\node [style=none] (473) at (6.3, -5.55) {\scriptsize{$i=9$}};
		\node [style=none] (474) at (2.975, -6.35) {\scriptsize{crossing\#3}};
		\node [style=none] (475) at (5.375, -6.55) {\scriptsize{$i=15$}};
		\node [style=reddot1] (476) at (-2.9, 0.75) {};
		\node [style=reddot1] (477) at (-4.1, -1.25) {};
		\node [style=none] (478) at (-3.25, 1.25) {};
		\node [style=none] (479) at (-3.25, 0.25) {};
		\node [style=none] (480) at (-4.75, -1) {};
		\node [style=none] (481) at (-4.25, -2.25) {};
		\node [style=none] (482) at (-1.25, -5.35) {\scriptsize{$e=-1$}};
		\node [style=none] (483) at (-0.575, -3.7) {\scriptsize{$e=1$}};
		\node [style=none] (484) at (-0.975, -2.55) {\scriptsize{$e=-1$}};
		\node [style=none] (485) at (-1, -1.325) {\scriptsize{$e=1$}};
		\node [style=none] (486) at (-1.35, -0.125) {\scriptsize{$e=1$}};
		\node [style=none] (487) at (-0.55, 1.9) {\scriptsize{$e=-1$}};
		\node [style=none] (488) at (-0.65, 3.25) {\scriptsize{$e=-1$}};
		\node [style=none] (489) at (-1, 4.5) {\scriptsize{$e=1$}};
		\node [style=none] (490) at (-1.875, -3.95) {};
		\node [style=none] (491) at (-1.525, -5.05) {};
		\node [style=none] (492) at (3.95, -3.025) {\scriptsize{$e=2$}};
		\node [style=none] (493) at (2.975, -0.95) {\scriptsize{$e=-2$}};
		\node [style=none] (494) at (2.95, 0.6) {\scriptsize{$e=2$}};
		\node [style=none] (495) at (3.25, 3.175) {\scriptsize{$e=-2$}};
		\node [style=none] (496) at (6.65, -5.075) {\scriptsize{$e=\pm0.2$}};
		\node [style=none] (497) at (5.65, -6.125) {\scriptsize{$e=\pm0.3$}};
		\node [style=none] (498) at (5.25, -7.15) {\scriptsize{$e=\pm0.1$}};
		\node [style=none] (499) at (-0.55, 0.35) {};
		\node [style=none] (500) at (1.05, -7.3) {};
		\node [style=none] (501) at (3.175, -3.2) {};
		\node [style=none] (502) at (2.325, -2.875) {};
	\end{pgfonlayer}
	\begin{pgfonlayer}{edgelayer}
		\draw (367.center) to (370.center);
		\draw [style=redline, in=15, out=-15, looseness=1.75] (416.center) to (414.center);
		\draw [style=redline, in=-180, out=-180, looseness=1.50] (411.center) to (415.center);
		\draw [style=redline, in=-150, out=0] (415.center) to (417.center);
		\draw [style=redline, in=15, out=150, looseness=1.75] (412.center) to (420.center);
		\draw (422.center) to (423.center);
		\draw [style=redline, in=60, out=30, looseness=1.75] (448.center) to (449.center);
		\draw [style=redline, in=150, out=-150] (448.center) to (456.center);
		\draw [style=redline, in=165, out=-150, looseness=1.25] (457.center) to (458.center);
		\draw [style=redline, in=-165, out=0, looseness=0.75] (458.center) to (459.center);
		\draw [style=redline, in=15, out=60, looseness=1.50] (460.center) to (459.center);
		\draw [style=redline, in=-135, out=15, looseness=0.75] (462.center) to (461.center);
		\draw [style=redline, in=45, out=-150, looseness=2.00] (463.center) to (464.center);
		\draw [style=redline, in=120, out=-15] (409.center) to (471.center);
		\draw [style=redline, in=-180, out=-60] (436.center) to (414.center);
		\draw [style=redline] (413.center) to (418.center);
		\draw [style=arrow2, bend right] (478.center) to (479.center);
		\draw [style=arrow2] (480.center) to (481.center);
		\draw [style=arrow2] (491.center) to (490.center);
		\draw [style=arrow2] (435.center) to (472.center);
		\draw [style=arrow2] (419.center) to (470.center);
		\draw [style=arrow2] (500.center) to (499.center);
		\draw [style=arrow2] (501.center) to (502.center);
	\end{pgfonlayer}
\end{tikzpicture}
        \caption{Multi-Component Crossing, Index, and Orientation Inputs}
        \label{fig14}
\end{wrapfigure}

The skein in Figure \ref{fig14} has two components ($l=2$). We set a reference point in the first component, assign an orientation direction (counterclockwise), and record the skein intersection information $E^{1}=[-1,1,-1,-0.2,0.3,-1,-1,1]$. We do the same for the second component, $E^{2}=[1,0.1,-2,2,-2,2,0.2,-0.3,-0.1,-1]$. Then, $E=[E^1,E^2]$. The index is given by $I^1=[1,2,3,9,15,4,7,8]$ for component 1, $I^2=[5,14,13,12,11,10,9,15,14,6]$ for component 2, and $I=[I^1,I^2]$. For orientation, $Q^{1}=[3,4,3,0,0,4,3,4]$, $Q^{2}=[3,0,4,5,4,5,0,0,0,4]$ and $Q=[Q^1]$.

%%%%%%%%%%%%%
\subsection{Skein Input Parameters}\label{skein_input}
%%%%%%%%%%%%%

In describing the initial skein, we produce three arrays: the crossing array $E_0$, the index array $I_0$, and the orientation array $Q_0$ as described above (NOTE: We are using the subscript 0 to indicate that this information is associated to the original skein, before the resolution of any crossings). We can include all this information in an initial array of the form: \\ \\
 $Z_0=[s_0,c_0,U_0,E_0,I_0,Q_0]$ \\ \\
\noindent where $s_0=0$ is the initial power of the variable $t$ (i.e., there is no variable $t$ prior to the application of the skein relations to resolve a crossing) and $c_0=1$ (the sign and value of the initial coefficient). The array $U_0$ is an array that stores the crossings to be resolved; as crossings are resolved elements are removed from $U$ until $U=\emptyset$. We thus recover the skein as $c_0t^{s_0}K_0 = 1t^0K_0$ where $K_0=[E_0,I_0,Q_0]$ is the complete description of the unique skein recovered by its crossing, index, and orientation arrays (as shown, up to a choice of orientation, in Proposition REF). At each stage in the skein computation, we resolve one of the $r$ crossings, producing $2^{r}$ terms. The final result, after the resolution of all crossings in the initial diagram, is a polynomial with $2^r$ diagrammatic terms, which can be described by array of the form:

\begin{equation*}
    Z=[Z_0,Z_1,...,Z_{2^r-1}]
\end{equation*}

where each component $Z_z$ represents one terms of the final Laurent polynomial, $Z_z= [s_z,c_z,U_z,E_z,I_z,Q_z]$. We recover a linear combination of skeins in the variable $t$ from $Z$ as:

\begin{equation*}
    c_0t^{s_0}K_0+c_1t^{s_1}K_1+...+c_{2^r-1}t^{s_{2^r-1}}K_{2^r-1}
\end{equation*}

where $K_z=[E_z,I_z,Q_Z]$, and the additional information ($c_z,s_Z$) stores the value of the coefficient in $\mathbb{C}[t,t^{-1}]$ as $c_zt^{s_z}$. Table 1 summarizes each of the entries in the skein array. The fact that the skein can be uniquely recovered from this the information given here is formalized in Proposition \ref{prop1}.

\begin{table}[H]
\label{Table_1}
\caption{Summary of the skein array}
\begin{center}
\begin{tabular}{|c|c|}
\hline
\multicolumn{2}{|c|}{$Z_z= [s_z,c_z,U_z,E_z,I_z,Q_z]$}\\
\hline
Skein Array Information: & Description: \\
\hline
$s_z$ & power of the variable $t$  \\ 
\hline
$c_z$ & sign and value of the constant, $s_zc_z \in \mathbb{C}[t,t^{-1}]$  \\ 
\hline
$U_z=[\pm1,\pm1,...,\pm1]$ & signs on crossings to be resolved ($dim(U_z)=r$)\\ 
\hline
$K_z=[E_z,I_z,Q_z]$ & skein input information\\
\hline
$E_z=[E_{z}^{0},E_{z}^{1},...E_{z}^{l},...,E_{z}^{k}]$ & skein crossing array  \\ 
\hline
$I_z=[I_{z}^{0},I_{z}^{1},...I_{z}^{l},...,I_{z}^{k}]$ & skein index array  \\ 
\hline
$Q_z=[Q_{z}^{1},Q_{z}^{2},...Q_{z}^{l}...,Q_{z}^{k}]$ & skein orientation array  \\ 
\hline
\end{tabular}
\end{center}
\end{table}

\begin{Prop} \label{prop1}
The resolution of crossings $1 \times 10^{-k_{r}}, 2 \times 10^{-k_{r}}, ..., w \times 10^{-k_{r}}$ by manipulations of the array $Z=[Z_0,Z_1,...,Z_z,...Z_{2^r-1}]$ produces the linear combination of skeins in $t$ with coefficients in $\mathbb{C}[t,t^{-1}]$:
\begin{equation}
    \sum_{z=0}^{2^w-1}c_zt^{s_z}K_z
\end{equation}

where $K_z=[E_z,I_z,Q_z]$ and $w$ is the integer corresponding to the label on the last crossing to be resolved. When $w=r$, 

\begin{equation}
    \sum_{z=0}^{2^r-1}c_zt^{s_z}K_z
\end{equation}

where $U=\emptyset$ and all crossings in the original crossings in the skein have been resolved.
\end{Prop}

Finally, we note that any skein embedded in a genus-2 handlebody can be uniquely recovered from the input information required in this paper ($E$, $I$, and $Q$). This is not the minimum representation of such a skein, but it is one which allows for automated computation and efficient manipulation.

\begin{Lem} \label{lem1}
    The representation $K=[E,I,Q]$ recovers a unique skein in a genus-2 handlebody. Any skein in a genus-2 handlebody has multiple representations, $K$, equivalent up to permutation of the array elements.
\end{Lem}

\proofname \; \;We recover the skein $K$ from the given arrays $E$ (skein intersection array), $I$ (index array), and $Q$ (orientation array). In each component $l$, we recover $K^{l}=[E^l,I^l,Q^l]$ as follows: the dimension of the arrays $dim(E^l)=dim(I^l)=dim(Q^l)=h$ is the number of points of interest we consider. Beginning with the array $I^l$ we label points along each strand of the base braid by their sequential index values, from bottom to top, by partitioning $I^l$ into three sets: $I^{l}_{1}$, $I^{l}_{2}$, and $I^{l}_{3}$ according to the corresponding values in $E^{l}$. $I^{l}_{1}$ contains all elements $i^{l}_{j}$ for which $|e^{l}_{j}|=1$. $I^{l}_{2}$ contains all elements $i^{l}_{j}$ for which $|e^{l}_{j}|=2$. $I^{l}_{3}=I^{l} \setminus (I^{l}_{1} \cup I^{l}_{2})$; thus $I^{l}_{3}$ contains all indices corresponding to skein self-intersections. This induces a partition on $E^{l}$ (into $E^{l}_{1}$, $E^{l}_{2}$, and $E^{l}_{3}$). Let $n_{1}=dim(I^{l}_{1})$, $n_{2}=dim(I^{l}_{2})$, and $n_{3}=dim(I^{l}_{3})$. The values of the elements of $I^{l}_{1}$ determine an ordering along strand 1. Similarly for $I^{l}_{2}$. Begin with the minimal element of $I^{l}_{1}$ (the element with the lowest index) and consider the corresponding value of $E^{l}_{l}$; WLOG, call this minimal element $i^{l}_0 \in I^{l}_1 \subset I^{l}$. It's corresponding value is $e^{l}_0 \in E^{l}_1 \subset E^{l}$. (Note: we assume that the array $E^{l}$ begins at the over/undercrossing of strand 1 with the lowest index; starting anywhere else, we can rearrange our array such that this element is the first entry by permuting the elements. The starting point we choose to begin recording skein information is arbitrary; we can always choose to record skein information beginning with the lowest index along strand 1, thus, we assume $|e^{l}_0|=1$). If $e^{l}_0=1$, we draw an overcrossing and if $e^{l}_0=-1$ we draw an undercrossing at this location. Returning to the original array $I^l$, beginning with $i^{l}_{0}$ we connect $i^{l}_{0}$ to $i^{l}_{1}$ by inspecting the value of $e^{l}_{1}$. Continuing in this fashion for all  $(e^{l}_{j},i^{l}_{j})$, we recover the unoriented skein. Now, fix any value $j$ and consider $(e^{l}_{j},i^{l}_{j},q^{l}_{j})$ and $(e^{l}_{j+1},i^{l}_{j+1},q^{l}_{j+1})$. The comparison of $q^{l}_{j}$ and $q^{l}_{j+1}$ determines the orientation of component $l$. Thus, the skein in uniquely determined for each component from the representation $K^{l}=[E^{l},I^{l},Q^{l}]$. \qedsymbol

\begin{Rem}
We note that each skein $K$ may have multiple equivalent representations. Indexing the skein based on height we only require the index to be sequential along each strand. The index labeling on the base braid does not affect the labeling on any other strand (independent between strands). Thus, the array $I^l$ used in each component $l$ is not unique.
\end{Rem}

%%%%%%%%%%%%
%%%%%%%%%%%%%%%%%%%
\section{Algorithms for Skein Computations}
%%%%%%%%%%%%%%%%%%%
%%%%%%%%%%%%

We first introduce an algorithm to resolve all crossings $1 \times 10^{-k_r}, 2 \times 10^{-k_r},..., r \times 10^{-k_r}$ in a given skein system. This is split between two sub-algorithms; the first (Algorithm 1) addresses the structure of the algorithm on array $Z$, resolving crossings one by one. The second (Algorithm 2) pulls out the details of that resolution on individual components of $Z$.

\begin{wrapfigure}{L}[0pt]{5.5cm}
    \centering
	\begin{tikzpicture}
	\begin{pgfonlayer}{nodelayer}
		\node [style=none] (11) at (0.25, 1.125) {};
		\node [style=none] (12) at (-1.5, 1.125) {};
		\node [style=none] (16) at (-1.5, -0.625) {};
		\node [style=none] (17) at (0.175, -0.55) {};
		\node [style=none] (23) at (-2, 0.25) {$t$};
		\node [style=none] (41) at (4.175, 1.35) {};
		\node [style=none] (42) at (4.15, -0.825) {};
		\node [style=none] (47) at (3.025, 1.35) {};
		\node [style=none] (48) at (3.025, -0.825) {};
		\node [style=none] (72) at (-4.825, -1.575) {$z_i$};
		\node [style=none] (74) at (-0.75, -1.625) {$z^{\prime}_{i}$};
		\node [style=none] (76) at (3.525, -1.65) {$z^{\prime\prime}_{i}$};
		\node [style=none] (78) at (2.475, 0.25) {$t^{-1}$};
		\node [style=none] (79) at (0.95, 0.25) {$+$};
		\node [style=none] (80) at (-3, 0.25) {$=$};
		\node [style=none] (81) at (-5.35, 1.175) {};
		\node [style=none] (82) at (-4.2, -0.825) {};
		\node [style=none] (83) at (-4.7, 0.3) {};
		\node [style=none] (84) at (-4, 1.25) {};
		\node [style=none] (85) at (-4.9, 0.125) {};
		\node [style=none] (86) at (-5.525, -0.825) {};
	\end{pgfonlayer}
	\begin{pgfonlayer}{edgelayer}
		\draw [bend left=90, looseness=1.50] (11.center) to (12.center);
		\draw [in=135, out=-135] (41.center) to (42.center);
		\draw [in=30, out=-30, looseness=0.75] (47.center) to (48.center);
		\draw [bend left=90, looseness=1.50] (16.center) to (17.center);
		\draw [in=90, out=-90, looseness=1.75] (81.center) to (82.center);
		\draw [in=-90, out=30, looseness=0.75] (83.center) to (84.center);
		\draw [in=75, out=-150] (85.center) to (86.center);
	\end{pgfonlayer}
\end{tikzpicture}
	\caption{$Z_i$ decomposes as a sum of $Z'_i$ and $Z"_i$, coefficients in $\mathbb{C}[t,t^{-1}]$}
        \label{fig15}
\end{wrapfigure}

 In Algorithms \ref{alg1} and \ref{alg2}, beginning with $Z=[Z_0]$ (the initial skein diagram, $1t^0K_0$ and continuing for $Z_i=c_it^{s_i}K_i$) the burden of our computation at each stage is to apply the skein relations to $Z_i$ given in Figure \ref{fig3} to produce two new skeins: $Z'_i$ and $Z''_i$. We then replace $c_it^{s_i}K_i$ by the combination $(c_i+c'_{i})t^{s_i+s'_{i}}K'_i+(c_i+c''_{i})t^{s_i+s''_{i}}K''_i$. The way we edit the array $Z_i$ to produce $Z'_i$ and $Z''_i$ (shown in Figure \ref{fig15}) depends on several factors: if the crossing to be resolved is a self-intersection in the same component $Z_i$ of the skein (Case 1) only that component needs to be manipulated. If, however, the crossing is between two components of the skein, $Z_i$ and $Z_j$, (Case 2) both components are manipulated and replaced by $Z'_i$ and $Z''_i$. 
 
 We perform these array manipulations (representing skein crossing resolution) via a series of operators: $\mathcal{C}$, $\mathcal{R}$, $\mathcal{F}$, and $\mathcal{B}$. In Algorithms 1 and 2 we resolve crossing $r$ in $E_i$ and then apply the same manipulations to $I_i$ and $Q_i$.  For ease of reference, we call the crossing to be resolved ``0.r'' (assuming $k_r=1$). In Definition \ref{def1}, we define these 4 operators, the input of which are arrays or pairs of arrays with given indices $j'$ and $j''$.

\begin{Def} \label{def1}
Given arrays $x$ ($dim(x)=m$), $y$ ($dim(y)=n$), $z$ ($dim(z)=k$), and indices $j'$, $j''$, we define the operators $\mathcal{C}$, $\mathcal{R}$, $\mathcal{F}$, $\mathcal{B}$ as follows:

\begin{center}
\begin{equation*}
    \mathcal{C}: x \xrightarrow[\text{split section, concatenate}]{\mathcal{C}} y, z 
\end{equation*}
\end{center}
\;\\
\begin{center}
$\mathcal{C}(x,j',j'') = y,z$ \\
$\text{where}\; y = [x_{j'+1},...,x_{j''-1}]$ 
$\text{and}\; z = [x_0,...,x_{j'-1},x_{j''+1}...,x_{m-1}]$, \\
\end{center}

\;\\
\begin{equation*}
    \mathcal{R}: x \xrightarrow[\text{reverse section}]{\mathcal{R}} z \\
\end{equation*}
\;\\
\begin{center}
    $\mathcal{R}(x,j',j'') = z$ \\
    $\text{where}\; z =[x_0,...,x_{j'-1},x_{j''-1},...,x_{j'+1},x_{j''+1}...,x_{m-1}]$, 
\end{center}

\;\\
\begin{equation*}
    \mathcal{F}: x,y \xrightarrow[\text{forward section insert}]{\mathcal{F}} z \\
\end{equation*}
\;\\
\begin{center}
    $\mathcal{F}(x,y,j',j'') = z$ \\
    $\text{where}\; z =[x_0,...,x_{j'-1},y_{j''+1},..., y_{n-1},y_0,...,y_{j''-1},x_{j'+1},...,x_{m-1}]$, 
\end{center}

\;\\
\begin{equation*}
    \mathcal{B}: x,y \xrightarrow[\text{backward section insert}]{\mathcal{B}} z \\
\end{equation*}
\;\\
\begin{center}
    $\mathcal{B}(x,y,j',j'') = z$ \\
    $\text{where}\; z =[x_0,...,x_{j'-1},y_{j''-1},...,y_0, y_{n-1},...,y_{j''+1},x_{j'+1},...,x_{m-1}]$. 
\end{center}

\end{Def}

 We adjust the array $U_i$ to record signs affected by the orientation changes in the recombination of skeins due to the resolution of a crossing. These changes affect future resolutions; the signs of crossings are only affected when the components in the array are reversed (i.e., when $\mathcal{R}$ and $\mathcal{B}$ are applied). As described in \ref{u}, we check for crossings in $E_i$ in the section of the skein that is reversed. In Case 1 ( i.e., crossings to be resolved are in the same component) we check for crossings between indices $j'$ and $j"$ (i.e., , any crossings in the array $[x_{j'+1},...,x_{j''-1}]$). If both the overcrossing and undercrossing of the same label lie in this section of the array, no change is needed. If only one part (either the overcrossing or the undercrossing) is in this section of the array we change the corresponding value in $U$ (specifically, the element of $U$ corresponding to sign of that crossing) to the opposite sign. In Case 2, the same condition holds except that we must check for crossings in the array $[y_0,...y_{n-1}]$.
 
\begin{algorithm}
	\caption{Resolve Crossings in Skein}
        \label{alg1}
	\begin{algorithmic}[1]
            \Statex Input: $Z$; $dim(Z)=1$
            \State $i=0$
            \While{$i<2^r$}
                \Comment{Resolve crossing $r$ in $Z_i$}
                \State $Z_i\xrightarrow{\text{Algorithm 2}} Z'_i, Z''_i$
                \State delete $Z_i$
                \State insert $Z'_{i}$ at index $i$
                \State insert $Z''_{i}$ at index $i+1$
                \State increase $i$ by 2
            \EndWhile
		\Statex Output: $Z$; $dim(Z)=2^r$
	\end{algorithmic}
\end{algorithm}

\begin{algorithm}
	\caption{Resolve Crossings in Skein}
        \label{alg2}
	\begin{algorithmic}[1]
            \Statex Input: $Z_i=[s_i,c_i,U_i,E_i,I_i,Q_i], r$
            \State Locate overcrossing $0.r$ and undercrossing $-0.r$ in $E_i$:
            \Comment{\emph{implement any array search for value}}
                \Statex $l' \leftarrow$ component of $\pm0.r$ in $E_i$
                \Statex $j' \leftarrow$ index of $\pm0.r$ in $E_i^{l'}$
                \Statex $l'' \leftarrow$ component of $\mp0.r$ in $E_i$
                \Statex $j'' \leftarrow$ index of $\mp0.r$ in $E_i^{l''}$
            \If{ $l'=l=l''$}
                \Comment{Case 1: over/undercrossing in same component}
                \State $EA_i: E_i^l \xrightarrow[\text{split section, concatenate}]{\mathcal{C}} \mathcal{C}_1(E_i^l),\mathcal{C}_2(E_i^l)$ 
                \State delete $E_i^l$ in $E_i$
                \State insert $\mathcal{C}_1(E_i^l)$ at index $l$ of $E_i$ 
                \State insert $\mathcal{C}_2(E_i^l)$ at index $l+1$ of $E_i$ 
                \State form corresponding $IA_i$, $QA_i$
                \Comment{Repeat Lines 3-6 for $IA_i$ and $QA_i$}
                \State adjust $U_i$ to form $UA_i$
                \Comment{Use Definition DEF}
                \State define $Z'_i=[s_i+1,c_i,UA_i,EA_i,IA_i,QA_i]$   
                \State $EB_i: E_i^l \xrightarrow[\text{reverse section}]{\mathcal{R}} \mathcal{R}(E_i^l)$
                \State delete $E_i^l$ in $E_i$
                \State insert $\mathcal{R}(E_i^l)$ at index $l$ of $E_i$
                \State form corresponding $IB_i$, $QB_i$
                \Comment{Repeat Lines 10-12 for $IB_i$ and $QB_i$}
                \State adjust $U_i$ to form $UB_i$
                \Comment{Use Definition DEF}
                \State define $Z''_i=[s_i-1,c_i,UB_i,EB_i,IB_i,QB_i]$  
            \Else
                \Comment{Case 2: over/undercrossing in different components}
                \State WLOG, assume $l'<l''$
                \State $EA_i: E_i^{l'},E_i^{l''} \xrightarrow[\text{forward section insert}]{\mathcal{F}} \mathcal{F}(E_i^{l'},E_i^{l''})$ 
                \State delete $E_i^{l''}$ in $E_i$
                \State delete $E_i^{l'}$ in $E_i$
                \State insert $\mathcal{F}(E_i^{l'}),E_i^{l''}$ at index $l'$ of $E_i$ 
                \State form corresponding $IA_i$, $QA_i$
                \Comment{Repeat Lines 18-21 for $IA_i$ and $QA_i$}
                \State adjust $U_i$ to form $UA_i$
                \Comment{Use Definition DEF}
                \State define $Z'_i=[s_i+1,c_i,UA_i,EA_i,IA_i,QA_i]$   
                \State $EB_i: E_i^{l'},E_i^{l''} \xrightarrow[\text{backward section insert}]{\mathcal{B}} \mathcal{B}(E_i^{l'},E_i^{l''})$ 
                \State delete $E_i^{l''}$ in $E_i$
                \State delete $E_i^{l'}$ in $E_i$
                \State insert $\mathcal{B}(E_i^{l'}),E_i^{l''}$ at index $l'$ of $E_i$ 
                \State form corresponding $IB_i$, $QB_i$
                \Comment{Repeat Lines 25-28 for $IB_i$ and $QB_i$}
                \State adjust $U_i$ to form $UB_i$
                \Comment{Use Definition DEF}
                \State define $Z''_i=[s_i+1,c_i,UB_i,EB_i,IB_i,QB_i]$ 
            \EndIf
		\Statex Output: $Z'_i$,$Z''_i$
	\end{algorithmic}
\end{algorithm}

\begin{algorithm}
	\caption{Induce Crossings in Skein}
        \label{alg3}
	\begin{algorithmic}[1]
            \Statex Input: $Z$; $dim(Z)=2^r$, where $r$ is the number of crossings in the original skein
            \State $k=0$
            \While{$k<dim(Z)$}
                \Statex Induce crossings in $Z_k$
                \Comment{$Z_k$ has no skein self-intersections, $U_k=\emptyset$ (output of Algorithm 1,2)}
                \For{$n=1,2$}
                    \State Form $P_n$, $N_n$, $B_n$ 
                    \Comment{Use Definition \ref{def2}}
                    \While{$P_n \neq \emptyset$}
                    \State $P_n, N_n, B_n \xrightarrow[\text{induction decision}]{\text{Algorithm 4}} (l',j'), (l'',j''), P'_n, N'_n, B'_n$
                    
                    \Comment{Use Algorithm 4 to locate where to induce crossing(s)}
                    
                    \Comment{$P_n$, $N_n$, $B_n$ are modified in Algorithm 4}
                    \State Update $P_n \leftarrow P'_n$, $N_n \leftarrow N'_n$, $B_n \leftarrow B'_n$
                    \State $Z_k, (l',j'), (l'',j'') \xrightarrow[\text{induce crossings}]{\text{Algorithm 5}} \hat{Z}_k$ 

                    \Comment{Use Algorithm 5 to induce crossings at location $(l',j')$ and $(l'',j'')$}

                    \Comment{\emph{$Z_k$ now has either one (0.1) or two (0.1, 0.2) crossings}}
                    \State Assign $r \leftarrow 1$ \Comment{Identify induced crossing $0.1$}
                    \State $\hat{Z}_k, r \xrightarrow[\text{resolve crossings}]{\text{Algorithm 2}} Z'_k, Z''_k$ 
                    \Comment{Use Algorithm 2 to resolve crossing 0.1}
                    \State delete $Z_k$
                    \State insert $Z'_k$ at index $k$
                    \State insert $Z''_k$ at index $k+1$
                    
                    \Comment{Lines 10-13: the kth term is replaced by the two new skein terms (resolved)}

                    \Comment{\emph{This increases $dim(Z)$ by 1; continue with the new $Z_k=Z'_k$}}
                    \If{x=2} 
                    \Comment{If 2 crossings were induced in Algorithm 5}
                    \State $r \leftarrow x$
                    \Comment{Identify induced crossing 0.2 in $Z'_K$ of line 12 \emph{(now at index k)}}
                    \State $Z_k, r \xrightarrow[\text{resolve crossings}]{\text{Algorithm 2}} Z'_k, Z''_k$ 
                    \Comment{Use Algorithm 2 to resolve crossing $0.2$}
                    \State delete $Z_k$ 
                    \State insert $Z'_k$ at index $k$
                    \State insert $Z''_k$ at index $k+1$
                    \Comment{\emph{This increases dim(Z) by 1}}
                    
                    \Comment{Identify induced crossing 0.2 in term $Z''_k$ of line 13 \emph{(now at index k+2)}}
                    \State $Z_{k+2}, r \xrightarrow[\text{resolve crossings}]{\text{Algorithm 2}} Z'_{k+2}, Z''_{k+2}$ 
                    \Comment{Use Algorithm 2 to resolve crossing $0.2$}
                    \State delete $Z_{k+2}$ 
                    \State insert $Z'_{k+2}$ at index $k+2$
                    \State insert $Z''_{k+2}$ at index $k+3$
                    
                    \Comment{Lines 14-23: Each of the terms from Lines 10-13 is replaced by two new terms (resolved)}
                    
                    \Comment{Thus, the 4 new terms are located at indices k, k+1, k+2, and k+3}
                    \EndIf
                    \Comment{Continue with modified $Z_k$ (and associated $P_n, N_n, B_n$) until $P_n=\emptyset$}
                    \EndWhile
                \EndFor
                \State $k+=1$
                \Comment{$Z_k$ is fully resolved; increase $k$ by 1 and consider next component ($Z_{k+1}$)}
            \EndWhile
		\Statex Output: $Z$; $dim(Z)>>2^r$, each $Z_k$ fully resolved as basis elements $x,y,z$ and the unknot
	\end{algorithmic}
\end{algorithm}
 
\begin{algorithm}
	\caption{Induction Decision}
        \label{alg4}
	\begin{algorithmic}[1]
            \Statex Input: $n$, $B_n$, $P_n$, $N_n$
            \State define $b=min(B_n)$, $c=min(P_n)$
            %\While{$P\neq\emptyset$}
                \While{$b=c$}
                    \Comment{Compare $min(P_n)$ to $min(B_n)$}
                    \State delete $c$ from $P_n$ $\rightarrow P'_n$
                    \Comment{$P'_n=P_n \setminus \{c\}$}
                    \State $c \leftarrow min(P_n)$ \Comment{recompute minimum of modified $P_n$}
                    \State delete $b$ from $B_n$ $\rightarrow B'_n$
                    \Comment{$B'_n=B_n \setminus \{b\}$}
                    \State $b \leftarrow min(B_n)$ 
                    \Comment{recompute minimum of modified $B_n$}
                \EndWhile
                \Comment{Since, $b \neq c$, $min(B_n) \neq min(P_n) \rightarrow min(B_n) = min(N_n)$}
                \State Define $M_n=N_n$
                \Comment{$M_n$ is a copy of $N_n$; $M_n$ is modified and discarded}
                \State $a=max(M_n)$ 
                \While{$a \geq c$}
                    \Comment{Compare $min(P_n)$ to $max(M_n)$}
                    \State delete $a$ from $M_n$
                    \State $a \leftarrow max(M_n)$ (recompute maximum of $M_n$)
                \EndWhile
                \Comment{now, $a < c$, $max(M_n)<min(P_n)$}
                
                \Comment{Induction Decision: switch indices $a=i_{j'}^{l'}$ and $c=i_{j''}^{l''}$}
            %\EndWhile
		\Statex Output: $j', l', j'', l''$ where $a=i_{j'}^{l'}$ and $c=i_{j''}^{l''}$, $P'_n, N_n, B'_n$
	\end{algorithmic}
\end{algorithm}

\begin{algorithm}
	\caption{Induce Crossings}
        \label{alg5}
	\begin{algorithmic}[1]
            \Statex Input: $Z_k$, $l', j', l'', j''$
            \If{$l'=l''=l \And |j''-j'|=1$}
                \Statex [same component (consecutive)]
                \If{$j'<j''$:}
                    \State insert $0.1$, $0$, $0$ at index $j''+1$ (after $j''$) of $E_k^l$, $I_k^l$, $Q_k^l$, respectively
                    \State insert $-0.1$, $0$, $0$ at index $j'$ (before $j'$) of $E_k^l$, $I_k^l$, $Q_k^l$, respectively
                    
                \Else
                    \State insert $-0.1$, $0$, $0$ at index $j'+1$ (after $j'$) of $E_k^l$, $I_k^l$, $Q_k^l$, respectively
                    \State insert $0.1$, $0$, $0$ at index $j''$ (before $j''$) of $E_k^l$, $I_k^l$, $Q_k^l$, respectively
                \EndIf
            \Else
                \Statex [different components or same component (non-consecutive)]
                \State insert $\pm0.1$, $0$, $0$ at index $j'+1$ of $E_k^{l'}$, $I_k^{l'}$, $Q_k^{l'}$, respectively
                \State insert $\mp0.2$, $0$, $0$ at index $j'+1$ of $E_k^{l'}$, $I_k^{l'}$, $Q_k^{l'}$, respectively
                \State insert $\mp0.1$, $0$, $0$ at index $j''+1$ of $E_k^{l'}$, $I_k^{l'}$, $Q_k^{l'}$, respectively
                \State insert $\pm0.2$, $0$, $0$ at index $j''$ of $E_k^{l''}$, $I_k^{l'}$, $Q_k^{l'}$, respectively
                \State Assign $x \leftarrow 2$ (Indicating 2 crossings have been introduced)
            \EndIf
            \State Update $u_k$ in $Z_k$ to reflect new signs:
		\Statex Output: $j', l', j'', l''$ where $a=i_{j'}^{l'}$ and $c=i_{j''}^{l''}$
	\end{algorithmic}
\end{algorithm}

For Algorithm \ref{alg3}, we implement an index sorting procedure along each strand ($n=1,2$). First we determine where to induce this crossing (``Induction Decision'' or Algorithm \ref{alg4}) by  forming three sub-arrays ($P$, $N$, and $B$) and comparing entries of $P$ and $N$ against $B$ to determine which indices to switch. We implement the following rule: overcrossings are moved to lower indices and undercrossings are moved to higher indices by either the application of an $R1$ or $R2$ Reidemeister move. This induces a crossing ($R1$) or pair of crossings ($R2$) at each index switch which is implemented in Algorithm \ref{alg5}. 

\begin{Def} \label{def2}
    The skein polynomial $Z$ is made up of individual skein terms $Z_k$, each containing a multi-curve $E_k$. Following the application of Algorithms 1-2, this multi-curve is made up of single curve components, $E_k^l$ with no-self intersections. Given a strand number $n$, we define unordered sets $P_n$, $N_n$, and $B_n$ from elements $(i_k^l)_j \in I_k^l$ based on the value of the corresponding elements $(e_k^l)_j \in E_k^l$ as follows:
    
    \begin{equation*}
        P_n=\{(i_k^l)_j\}, \text{such that } (e_k^l)_j=n,
    \end{equation*}
    
    \begin{equation*}
        N_n=\{(i_k^l)_j\}, \text{such that } (e_k^l)_j=-n,
    \end{equation*}

    \begin{equation*}
        B_n=\{(i_k^l)_j\}, \text{such that } |(e_k^l)_j|=n.
    \end{equation*}
\end{Def}

This amounts to partitioning the indices of the multi-curve $E_k$ into indices with a positive crossing on strand $n$ ($P_n$), indices with a negative crossing on strand $n$ ($N_n$), and a set that contains all indices along strand $n$ ($B_n$). Note that this partition is a single set across all single curve components $E_k^l$ (rather than a separate set for each component $l$). 

The induced crossings are resolved by application of Algorithms \ref{alg1} and \ref{alg2}. The result of Algorithms \ref{alg3} and \ref{alg4} is the resolution of the skein into a linear combination of the fundamental curves given in REF. This is proved in Theorem \ref{thm1} below.

\begin{Thm} \label{thm1}
    The strand sorting method of Algorithm \ref{alg3} produces only basis elements of the Kauffman Bracket Skein Module of the genus-2 handlebody: $x,y,z$ and the unknot.
\end{Thm}

\proofname \;The same method is applied to both strands of the base braid separately so consider elements along the first strand of the base braid ($|e^l_j|=1$). The result of the strand sorting method is that all overcrossings are moved to lower indices and undercrossings are moved to higher indices. This is equivalent to requiring the following condition along each strand: each skein term has the property that for all entries $e^l_j=1$ and all entries $e^{l'}_{j'}=-1$ the corresponding indices satisfy the relation $i^l_j<i^{l'}_{j'}$. We then realize that a skein with this requirement in the genus 2 handlebody is equivalent to a skein in the 2-punctured disk and all non-intersecting skeins in the 2-punctured disk are either $x$, $y$, $z$, or the unknot. \qedsymbol

Consecutive index elements along a strand are considered pairwise from bottom to top. If the pair is a pair of overcrossings (or a pair of undercrossings), no change is made. If the pair is one overcrossing and one undercrossing, and the overcrossing has a higher index, a crossing is induced by a Reidemeister $R2$ move between the two strands, switching their positions and thus, their indices. If the overcrossing and undercrossing pair are in the same component, a single crossing is induced by a Reidemeister $R1$ move. The crossing(s) are resolved, producing either 2 or 4 skein terms in which the indices of the pair have been switched (the overcrossing now has a lower index than the undercrossing) or the same pair no longer exists (if the skein has been reduced). Continuing in this manner, all such pairs are switched until the condition above is satisfied. 

%%%%%%%%%%%%%%
%%%%%%%%%%%%%%%%%%%%%%
\section{Conclusion and Future Work}
%%%%%%%%%%%%%%%%%%%%%%
%%%%%%%%%%%%%%

We produced a method to resolve all crossings in a skein embedded in a genus-2 handlebody (Algorithm \ref{alg1} and \ref{alg2}) and then methodically induce crossings on the resulting linear combination of skeins (Algorithms \ref{alg3} and \ref{alg4}) to produce Laurent polynomial in the Pryzycki basis $x,y,z$ of the Kauffman Bracket Skein Module. This method, called "Resolve Original Diagram; Induce Crossings After (RODICA)" is named in memory of Dr. Rodica Gelca (Texas Tech University, Lubbock TX). This method provides a robust way to encode diagrammatic skein information and a means for the automation of all skein computations in a genus-2 handlebody. We further developed this method for the automation of skein computations in the complement of all 2-bridge knots, the topic of a follow on paper "Algorithms for Skein Manipulation and Automation of Skein Computations."

\clearpage
%%%%%%%%%%%%%%%%%%%%%%%%%%%%%%%%%%%%%%%%%%%%
%2. Using natbib package (you may use other bibliography types). Refer to http://en.wikibooks.org/wiki/LaTeX/Bibliography_Management
\bibliographystyle{ieeetr}    %% this one looks best
%\bibliographystyle{apalike}   %% looks okay for dissertations but it puts quotes around titles in references
%\setcitestyle{authoryear, open={((},close={))}}
\bibliography{SkeinAlgBib} %SampleBib

\begin{thebibliography}{1}

\bibitem{jones1985}
{Vaughan F. R. Jones}, ``{A polynomial invariant for knots via von Neumann
  algebras},'' {\em Bulletin (New Series) of the American Mathematical
  Society}, vol.~12, no.~1, pp.~103 -- 111, 1985.

\bibitem{kauffman1987}
{Louis H. Kauffman}, ``{State models and the Jones polynomial},'' {\em
  Topology}, vol.~26, no.~3, pp.~395--407, 1987.

\bibitem{witten1989}
{Edward Witten}, ``Quantum field theory and the jones polynomial,'' {\em
  Communications in Mathematical Physics}, vol.~121, no.~3, pp.~351--399, 1989.

\bibitem{reshetikhin1991}
{Nicolai Reshetikhin} and {Vladmir G Turaev}, ``{Invariants of 3-manifolds via
  link polynomials and quantum groups},'' {\em Inventiones mathematicae},
  vol.~103, no.~3, pp.~547--598, 1991.

\bibitem{turaev2016}
{Vladimir G. Turaev}, {\em {Quantum Invariants of Knots and 3-Manifolds}}.
\newblock Berlin, Boston: De Gruyter, 2016.

\bibitem{blanchet1992}
{Christian Blanchet}, {Nathan Habegger}, {Gregor Masbaum}, and {Pierre Vogel},
  ``Topological quantum field theories derived from the kauffman bracket,''
  {\em Topology}, vol.~31, pp.~685--699, 1992.

\bibitem{lickorish1993}
{W. B. Raymond Lickorish}, ``{The skein method for three-manifold
  invariants},'' {\em Journal of Knot Theory and Its Ramifications},
  vol.~2(02), pp.~171--194, 1993.

\bibitem{przytycki1991}
{Jozef H Przytycki}, ``{Skein modules of 3-manifolds},'' {\em Bulletin of the
  Polish Academy of Sciences}, pp.~91--100, 1991.

\end{thebibliography}

\end{document}